\documentclass[11pt,a4paper,reqno]{amsart}
\usepackage{amsmath,amsthm}
\usepackage{amssymb,latexsym}
\usepackage{enumerate}
\usepackage{graphicx}
\usepackage{subcaption}
\newtheorem{theorem}{Theorem}[section]

\newtheorem{example}[theorem]{Example}

\newtheorem{remark}[theorem]{Remark}

\makeatletter
\let\mytagform@=\tagform@
\def\tagform@#1{\maketag@@@{\bfseries(\ignorespaces#1\unskip\@@italiccorr)}\hspace{-40pt}}
\renewcommand{\eqref}[1]{\textup{\mytagform@{\ref{#1}}}}
\makeatother

\begin{document}
\title{\bf A queueing/inventory and an insurance risk model}
\author{Onno Boxma, Rim Essifi and Augustus J.E.M. Janssen}
\thanks{Department of
Mathematics and Computer Science, Eindhoven University of Technology}

\maketitle

\date{\today}

\begin{abstract}
We study an $M/G/1$-type queueing model with the following additional feature.
The server works continuously, at fixed speed, even if there are no service requirements.
In the latter case, it is building up inventory, which can be interpreted as negative workload.
At random times, with an intensity $\omega(x)$ when the inventory is at level $x>0$,
the present inventory is removed, instantaneously reducing the inventory to zero.
We study the steady-state distribution of the (positive and negative) workload levels
for the cases $\omega(x)$ is constant and $\omega(x) = a x$.
The key tool is the Wiener-Hopf factorisation technique. When $\omega(x)$ is constant, no specific assumptions will be made on the
service requirement distribution. However, in the linear case, we need some algebraic hypotheses concerning the
Laplace-Stieltjes transform of the service requirement distribution.
Throughout the paper, we also study a closely related model coming from insurance risk theory.\\
\ \\
Keywords: M/G/1 queue, Cram\'er-Lundberg insurance risk model, workload, inventory, ruin probability, Wiener-Hopf technique.\\
\ \\
2010 Mathematics Subject Classification: 60K25, 90B22, 91B30, 47A68.
\end{abstract}

\section{Introduction}
In this paper, we study two related stochastic models: a queueing/inventory model and an insurance risk model.
The insurance risk model is a relaxation of the classical Cram\'er-Lundberg model.
Unlike that model, when the capital of the insurance company becomes negative, the company continues to operate in the same way.
However, during periods of negative surplus, the company can go {\em bankrupt}. It goes bankrupt according to some bankruptcy rate $\omega(x)$
when the negative surplus equals $x$; see Figure 1.
\begin{figure}
\begin{center}
\includegraphics[width=0.7\textwidth]{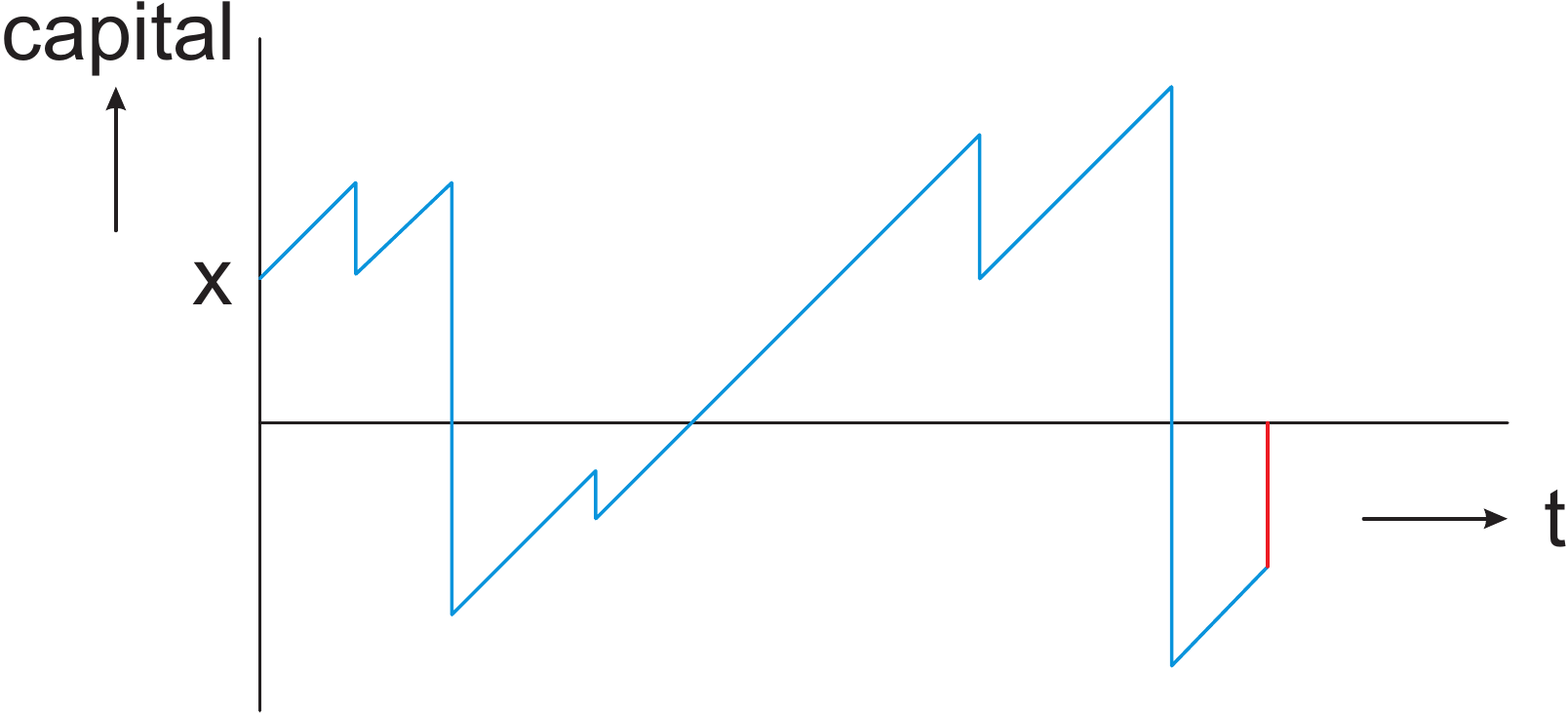}
\caption{The capital of an insurance company}
\end{center}
\end{figure}
This relaxation of the ruin concept was introduced in \cite{AlbGS}, and studied in \cite{AlbL}
for exponential claim sizes and various bankruptcy rates.
One of our goals in the present paper is to extend some of the results in \cite{AlbL}
to {\em{general claim size distributions}}. In particular, we aim to study the bankruptcy probability
when starting at $x$, for both positive {\em and} negative values of $x$.

It is well-known that the Cram\'er-Lundberg model is dual to the $M/G/1$ queueing model
with the same arrival rate and with a service time distribution that equals the claim size distribution in the Cram\'er-Lundberg model.
More precisely, cf.\ \cite{AA}, p. 52: The probability of ruin in the Cram\'er-Lundberg model with initial capital $x$
equals the probability that the steady-state virtual waiting time (or workload) in the $M/G/1$ queue exceeds $x$.
This has led us to think about queueing models that are relaxations of the $M/G/1$ queue in a similar way as
the Albrecher-Lautscham bankruptcy model is a relaxation of the Cram\'er-Lundberg model; see also \cite{AlbrecherBoxmaKuijstermans}.

\indent In this paper we study  the following queueing/inventory model.
Customers arrive according to a Poisson process, and require i.i.d.\ service times.
When there are customers, the server works at unit speed.
So far, this is the $M/G/1$ setting.
However, when there are no customers, the server still keeps on working at the same speed.
In that way, it is building up inventory. During periods in which there are no customers, inventory is instantaneously
removed according to a Poisson process with rate $\omega(x)$ when the amount of inventory equals $x$.
That inventory is, e.g., sold. The server just keeps on working; and when a customer arrives and
its service request can be satisfied from the inventory, then that is done instantaneously.
See Figure 2. In \cite{AlbrecherBoxmaKuijstermans} this two sided queueing/inventory model has been analysed for the case of exponentially
distributed service times.
Our queueing/inventory model is related to classical $M/G/1$ and inventory models.
An important inventory model is the basestock model, in which a server produces products
until the inventory has reached a certain basestock level, with requests for products arriving according to a Poisson process.
A request that cannot immediately be satisfied joins a backorder queue.
However, that model has a finite basestock level,
and hence essentially differs from our model.
Two papers which are to some extent related to our paper are \cite{BPP} and \cite{PSZ}.
These consider a production/inventory model with a so-called sporadic clearing policy.
The system is continuously filled at fixed rate, and satisfies demands at Poisson epochs. Under the sporadic clearing policy,
clearing of all inventory takes place at a random time (which in \cite{PSZ} is independent of the content process).
The authors obtain explicit results for an expected discounted cost functional.
\cite{PSZ} allows the demands to have a general distribution, whereas these are exponentially distributed in \cite{BPP}.
\begin{figure}
\begin{center}
\includegraphics[width=0.7\textwidth]{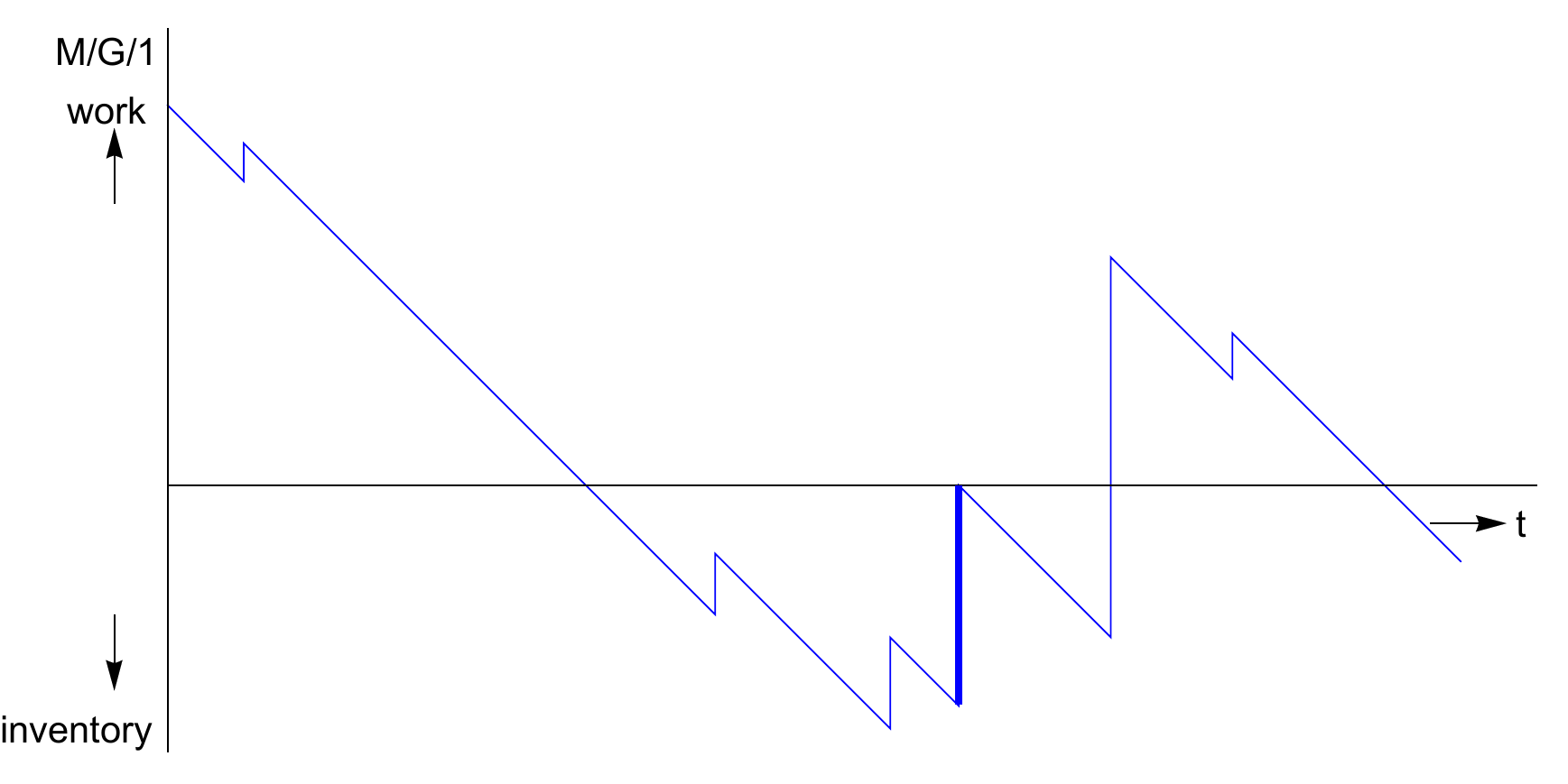}
\caption{Work and inventory in a queueing/inventory model}
\end{center}
\end{figure}
The main contributions of our paper are (i) an exact analysis of the Albrecher-Lautscham bankruptcy model
and the queueing/inventory model,
with generally distributed claim sizes, respectively generally distributed service requirements,
for the case of $\omega(x) \equiv \omega$; and (ii) a detailed analysis of the queueing/inventory model
for the case of $\omega(x) = ax$. The latter case turns out to lead to an inhomogeneous first order differential equation
with removable singularities, and its analysis gives rise to intricate calculations that are given in an appendix.
A key tool that we are using in the analysis of the two models is Wiener-Hopf factorization.
The results of the present paper might be used for optimization purposes; for example, one might try
to choose $\omega$ or $a$ (in the case $\omega(x) = ax$) such that a particular objective function is optimized.

The paper is organized as follows.
The queueing/inventory and insurance risk model are both described in detail in Section~\ref{modeldes}.
Integral equations for the main performance measures (workload and inventory densities in the queueing/inventory model,
bankruptcy probability when starting at level $x$ in the insurance risk model) are presented in Section~\ref{maineq}.
In Section~\ref{analysis} these equations are solved for $\omega(x) \equiv \omega$ and general service requirement distribution, respectively
general claim size distribution.
The queueing/inventory model is treated in Section~\ref{linearomega} for the case $\omega(x) = ax$.
That analysis makes use of Laplace transforms and complex analysis.
We assume that the service requirements having a rational Laplace transform; see also the appendix, that is devoted
to a detailed exposition of technical complications which arise when one tries to determine the Laplace transform
of the density of the inventory level by solving a first order algebraic differential equation with singularities.
Section~\ref{direct} considers that very same case, under the assumption of exponentially distributed service requirements,
without resorting to Laplace transforms.

\section{Model description}
\label{modeldes}
\subsection{Queueing/inventory model}
\label{Qmodel}
We study the following model (cf.\ Figure 2). Customers arrive according to a Poisson process with rate $\lambda$. Their service requirements
are independent, identically distributed random variables $B_1,B_2,\ldots$ with common distribution $B(.)$ and LST (Laplace-Stieltjes Transform)
$\beta(\cdot)$. The server works continuously, at a fixed speed which is normalized to $1$ - even if there are no service requirements. In the latter case, the server is building up inventory, which can be interpreted as negative workload. At random
times, with an intensity $\omega(x)$ when the inventory is at level $x>0$, the present inventory is removed, instantaneously reducing the inventory to zero. Put differently, inventory is removed according to a Poisson process with a rate that depends on the amount of inventory present.\\
Denote the required work per time unit by $\rho:=\lambda \mathbb EB$. We assume that $\rho<1$. This ensures that the steady state workload distribution exists. Let $V_+(x)$, $x>0$ denote this steady-state workload distribution, and $v_+(x)$ its density.\\
During the times in which the inventory level is positive, there is an upward drift $1-\rho$ of that inventory level; but when $\omega(x)>0$
for $x$ sufficiently large, the inventory level will always eventually return to zero, and the steady-state inventory distribution will exist.
Let $V_-(x)$, $x>0$ denote this steady state distribution, and $v_-(x)$ its density.\\

\subsection{Insurance risk model}
\label{IRmodel}
The problem that we will deal with here was introduced by Albrecher, Gerber and Shiu \cite{AlbGS} and investigated by Albrecher and Lautscham in \cite{AlbL}.
We will examine the bankruptcy probability for a surplus process with jumps.
Consider a Cram\'er-Lundberg setup to describe the insurer's surplus $C_t$ at time $t$ as
\begin{equation}\label{Cramer_lindberg_surplus}
C_t=x+ct-S_t,
\end{equation}
where $C_0=x$ is the initial surplus, $c$ is the premium rate, and $S_t$ is the aggregate claim amount up to time $t$ modeled as a compound Poisson
process with intensity $\lambda$ and positive jump sizes $Y_1,Y_2,\dots$ with cumulative distribution function $F_Y(\cdot)$.
In order to compare the results for the bankruptcy model with those for the queueing/inventory model,
we shall take $F_Y(\cdot) = B(\cdot)$, so its LST is $\beta(s)$.
It is assumed here that the insurer may be allowed to continue the business despite a temporary negative surplus.
Concretely, consider a suitable locally bounded bankruptcy rate function $\omega(-C_s)$ depending on the size of the negative surplus $C_s<0$.
If no bankruptcy event has occurred yet at time $s$, then the probability of bankruptcy in the time interval $[s,s+dt)$ is $\omega(-C_s)dt$.
We assume that $\omega(\cdot)\geq0$ and $\omega(x)\geq \omega(y)$ for $|x|\geq |y|$ to reflect that the likelihood of bankruptcy does not decrease as the surplus becomes more negative.
Let $\tau$ be the resulting time of bankruptcy, and define the overall probability of bankruptcy as
\begin{equation}\label{Def_u}
u(x)=\mathbb E[1_{\{\tau<\infty\}}\vert C_0=x]=\mathbb P[\tau <\infty\vert C_0=x].
\end{equation}
The idea is that whenever the surplus level becomes negative, there may still be a chance to survive, and  survival is less likely the lower such a negative surplus is.
For $x>0$, set $u_+(x):=u(x)$, $\tilde{u}_-(x):=u(-x)$ and $u_-(x):=1-\tilde{u}_-(x)$.\\

\section{Main equations}
\label{maineq}
In this section we present integral equations for the main performance measures (workload and inventory densities in the queueing/inventory model,
bankruptcy probability when starting at level $x$ in the insurance risk model).
\subsection{Queueing model}\ \\
The level crossing technique \cite{Brill}
yields the following integral
equations for the workload and inventory densities, by equating the rates
at which level $x$ is downcrossed and upcrossed in steady state:
\begin{equation}\label{queuing model>0}
v_+(x)=\lambda \int_0^{x} \mathbb
P(B>x-y)v_+(y)dy+\lambda\int_0^{+\infty}\mathbb P(B>x+y)v_-(y)dy,\ x>0,
\end{equation}
\begin{equation}\label{queuing model<0}
v_-(x)=\lambda \int_x^{+\infty} \mathbb
P(B>y-x)v_-(y)dy+\int_x^{+\infty}\omega(y)v_-(y)dy,\ x>0.
\end{equation}
We introduce the Laplace transforms
$\phi_+(s):=\int_0^{+\infty}e^{-sx}v_+(x)dx$ and
$\phi_-(s):=\int_0^{+\infty}e^{-sx}v_-(x)dx$ for $\rm {Re}\  s
\geq 0$. Multiplying both sides of
(\ref{queuing model>0}) with $e^{-sx}$ for $\rm {Re}\  s
\geq 0$ and both sides of
({\ref{queuing model<0}}) with $e^{sx}$ for $\rm {Re}\  s
\leq 0$, integrating and adding
both equations gives, after some calculations:
\begin{equation}{\label{eq_Wiener_Hopf_0}}
[1-\lambda \frac{1-\beta(s)}{s}][\phi_+(s)+\phi_-(-s)]=1/s
\int_{0}^{+\infty}(e^{sy}-1)\omega(y)v_-(y)dy, \quad \text{for}
\quad \rm {Re}\ s=0.
\end{equation}
\subsection{Insurance model}\ \\
According to \cite{AlbL}, one can write
\begin{eqnarray}
0&=&c  u'_+(x)-\lambda u_+(x)+\lambda \left(\int_0^x u_+(x-y)dB(y)+\int_x^{+\infty}u_-(y-x)dB(y)\right),\quad x>0, \label{Insurance_model_equation_+}\\
0&=&-c \tilde{u}_-'(x)-(\lambda+\omega(-x)) \tilde{u}_-(x)+\omega(-x)+\lambda \int_0^{+\infty}\tilde{u}_-(x+y)dB(y),\quad x>0. \label{Insurance_model_equation_-}
\end{eqnarray}
 Adding and subtracting $1$, Equation (\ref{Insurance_model_equation_-}) is equivalent to
\begin{equation}\label{Insurance_model_equation_-_bis}
0=-c  u_-'(x)+\lambda u_-(x)+\omega(-x)u_-(x)-\lambda \int_0^{+\infty}u_-(x+y)dB(y), \quad x>0.
\end{equation}
One can dominate the function $u_+$
by the classic ruin function and under assumptions on the existence of the second moment of $B(\cdot)$, the classic ruin function
is integrable. Hence, the function $u_+$ is integrable. Similarly, one can argue that the function $u_-$ is integrable.
For $\rm {Re}\  s\geq 0$, introduce the Laplace transforms
$\Psi_+(s):=\int_0^{+\infty}e^{-sx}u_+(x)dx$,
$\Psi_-(s):=\int_0^{+\infty}e^{-sx}u_-(x)dx$ and  $\beta(s):=\int_0^{+\infty}e^{-sy}dB(y)$.\\
Multiply both sides of (\ref{Insurance_model_equation_+})
with $e^{-sx}$ for $\rm {Re}\  s\geq 0$ and both sides of (\ref{Insurance_model_equation_-_bis})
with $e^{sx}$ for $\rm {Re}\  s\leq 0$; integrate and add
both equations for ${\rm{Re}}\ s=0$. After some calculations and using the fact that the continuity
of the function $u$ in zero implies $u_+(0)=1-u_-(0)$, we obtain:
\begin{eqnarray}
&\ &\left(\lambda \beta(s)+cs-\lambda\right)\Psi_+(s)+\frac{\lambda}{s}(1-\beta(s))-c \label{Wiener-Hopf_insurance_1}\\
&=&\left(\lambda \beta(s)+cs-\lambda\right)\Psi_-(-s)-\int_0^{+\infty}\omega(-x)u_-(x)e^{sx}dx, ~~~{\rm Re}~s = 0. \nonumber
\end{eqnarray}
In the next section, we restrict ourselves to the case where the
function $\omega(\cdot)$ is constant.

\section{Analysis for $\omega(\cdot)$ constant}
\label{analysis}
\subsection{Queueing model}\ \\
Assume in this section that the function $\omega(\cdot)$ introduced in Subsection~\ref{Qmodel} is constant, i.e., there exist $\omega>0$ such that
for all $x\geq 0$, one has $\omega(x)=\omega$.
Equation ({\ref{eq_Wiener_Hopf_0}}) becomes
\begin{equation}{\label{eq_Wiener_Hopf_1}}
\left[s-\lambda
(1-\beta(s))\right]\phi_+(s)=\left[\omega-\left(s-\lambda(1-\beta(s))\right)\right]\phi_-(-s)-\omega\phi_-(0)
 , \quad \text{for}
\quad \rm {Re}\ s=0.
\end{equation}
We are going to determine both unknown functions $\phi_+(s)$ and $\phi_-(-s)$ for $\rm {Re}\ s\geq0$ by formulating and solving
a Wiener-Hopf problem (cf.\ \cite{CohenWH}). A key step in this procedure is to rewrite Equation (\ref{eq_Wiener_Hopf_1}) such that the lefthand side is
analytic on $\rm {Re}\ s>0$ and the righthand side is analytic on $\rm {Re}\ s<0$. Liouville's theorem can subsequently be used
to identify the lefthand and righthand side.\\
Set
$f_{\omega,\lambda}: s \mapsto \lambda \beta(s) +s-\lambda-\omega, s\geq 0$ and $f_{0,\lambda}: s\mapsto \lambda \beta(s) +s-\lambda, s\geq 0$.
According to \cite{Cohen} p.\ $548$, the constant $\omega$ being positive, the function $f_{\omega,\lambda}$ has only one zero
$s=\delta(\omega,\lambda)$ and this zero is simple satisfying $\rm {Re}\ \delta(\omega,\lambda)>0$.
In fact, as $\omega$ is real in our case, a plot of $\omega + \lambda(1-\beta(s))$ versus $s$ immediately shows that this zero $\delta(\omega,\lambda)$ is real.
Also, the function $f_{0,\lambda}$ has $s=0$ as its only zero and this zero is simple. In particular the functions
$g_{\omega,\lambda}:\quad s \mapsto \frac {f_{\omega,\lambda}(s)}{s-\delta(\omega,\lambda)}$
and $g_{0,\lambda}: s \mapsto \frac {f_{0,\lambda}(s)}{s}$ are analytic on ${\rm{Re}}\ s>0$, continuous on ${\rm{Re}}\ s\geq0$ and take non-zero values on ${\rm{Re}}\ s\geq0$.\\
One can rewrite Equation ({\ref{eq_Wiener_Hopf_1}}) for $\rm {Re}\ s=0$ as:
 \begin{equation}{\label{eq_Wiener_Hopf_2}}
\frac{\left(s-\delta(\omega,\lambda)\right)\left[s-\lambda
(1-\beta(s))\right]\phi_+(s)}{s-\lambda(1-\beta(s))-\omega}+
\frac{\omega\phi_-(0)\left(s-\delta(\omega,\lambda)\right)}{s-\lambda(1-\beta(s))-\omega}
= -\left(s-\delta(\omega,\lambda)\right)\phi_-(-s).
\end{equation}
We now use the Wiener-Hopf factorization technique. The lefthand side of ({\ref{eq_Wiener_Hopf_2}}) is analytic
on $\rm {Re}\ s>0$ and continuous on $\rm {Re}\ s\geq0$; the righthand side is analytic
on $\rm {Re}\ s<0$ and continuous on $\rm {Re}\ s\leq0$. In addition, both sides coincide on $\rm {Re}\ s=0$. Then by Liouville's theorem
(cf.\ \cite{Titchmarsh}, p.$85$), there exist $ n \geq 0$
and a polynomial $R_n(s)$ of degree $n$ such that:
 \begin{eqnarray}{\label{factors_identification}}
 -\left(s-\delta(\omega,\lambda)\right)\phi_-(-s)&=&R_n(s)\ \text{for}\ \rm {Re}\ s\leq0 {\label{factors_identification_1}},\\
 \frac{\left(s-\delta(\omega,\lambda)\right)\left[s-\lambda
(1-\beta(s))\right]\phi_+(s)}{s-\lambda(1-\beta(s))-\omega}+
\frac{\omega\phi_-(0)\left(s-\delta(\omega,\lambda)\right)}{s-\lambda(1-\beta(s))-\omega}&=&R_n(s)\ \text{for}\ \rm {Re}\ s\geq0. {\label{factors_identification_2}}
 \end{eqnarray}
 Using Equation ({\ref{factors_identification_1}}) and the fact that $\displaystyle \lim_{s\rightarrow -\infty}\phi_-(-s)=0$, one has ${\rm deg} (R_n(s))=0$; say $R_n(s)=A$ where $A \in \mathbb C$. One gets
  \begin{equation}{\label{identification_phi-(-s)1}}
 \phi_-(-s)=\frac{A}{\delta(\omega,\lambda)-s}\ \text{for}\ \rm {Re}\ s\leq0,
 \end{equation}
 in particular
 \begin{equation}{\label{identification_phi-(0)}}
 \phi_-(0)=\frac{A}{\delta(\omega,\lambda)}.
 \end{equation}
 On the other hand, Equation ({\ref{factors_identification_2}}) yields
 \begin{equation}{\label{factors_identification_ph_+}}
 \phi_+(s)=\frac{A\left(s-\lambda(1-\beta(s))-\omega\right)}
 {\left(s-\lambda(1-\beta(s))\right)\left(s-\delta(\omega,\lambda)\right)}
-\frac{\omega\phi_-(0)}{s-\lambda\left(1-\beta(s)\right)}.
\end{equation}
After some calculations and using the notations introduced above and Equation ({\ref{identification_phi-(0)}}), one gets
\begin{equation}\label{identification_Phi_+_interm}
 \phi_+(s)=\frac{A\left(g_{\omega,\lambda}(s)-\frac{\omega}{\delta(\omega,\lambda)}\right)}{sg_{0,\lambda}(s)}.
 \end{equation}
 We now calculate the unknown constant $A$, and through Relations ({\ref{identification_phi-(0)}}) and (\ref{identification_Phi_+_interm}) we determine the functions $\phi_+$ and $\phi_-$.
 Note that $g_{\omega,\lambda}(0) =\frac{-f_{\omega,\lambda}(0)}{\delta(\omega,\lambda)}=\frac{\omega}{\delta(\omega,\lambda)}$, therefore one can write for $ \rm {Re}\ s\geq0$:
 \begin{equation}{\label{factors_identification_ph_+_1}}
 \phi_+(s)=\frac{A}{g_{0,\lambda}(s)}\frac{g_{\omega,\lambda}(s)-g_{\omega,\lambda}(0)}{s},
 \end{equation}
 the function $s \mapsto \frac{g_{\omega,\lambda}(s)-g_{\omega,\lambda}(0)}{s}$ clearly being analytic for $ \rm {Re}\ s>0$ and continuous for
 $ \rm {Re}\ s\geq0$.
\\
One has
 \begin{equation}
\phi_+(0)=\lim_{s\rightarrow 0}\phi_+(s)=\lim_{s\rightarrow 0} \frac{A}{g_{0,\lambda}(s)}\frac{g_{\omega,\lambda}(s)-g_{\omega,\lambda}(0)}{s}.
 \end{equation}
But
\begin{equation}
g_{0,\lambda}(0)=\lim_{s\rightarrow 0}\frac{f_{0,\lambda}(s)}{s}
 =\lim_{s\rightarrow 0} 1+ \lambda \frac{\beta(s)-1}{s}=1-\lambda\mathbb E(B)=1-\rho.
\end{equation}
One also has
$\displaystyle \lim_{s\rightarrow 0} \frac{g_{\omega,\lambda}(s)-g_{\omega,\lambda}(0)}{s}=g'_{\omega,\lambda}(0)
=\frac{\omega-(1-\rho)\delta(\omega,\lambda)}{\delta^2(\omega,\lambda)}$.
Then,
\begin{equation}{\label{factors_identification_ph_+_2}}
\phi_+(o)= A \frac{\omega-(1-\rho)\delta(\omega,\lambda)}{(1-\rho)\delta^2(\omega,\lambda)}.
\end{equation}
Using the relation $ \phi_-(0)+\phi_+(0)=1$, and Equations ({\ref{identification_phi-(0)}})
and ({\ref{factors_identification_ph_+_2}}), one gets
$$ A [ \frac{1}{\delta(\omega,\lambda)}+\frac{\omega-(1-\rho)\delta(\omega,\lambda)}{(1-\rho)\delta^2(\omega,\lambda)}]=1,$$
 which implies
 \begin{equation}{\label{Calcul_A}}
 A=\frac{(1-\rho)\delta^2(\omega,\lambda)}{\omega}.
 \end{equation}
 Finally we obtain the following expressions for $\phi_-(s)$ for $\rm {Re}\ s\leq0$ and for
 $\phi_+(s)$ for $\rm {Re}\ s\geq0$:
 \begin{eqnarray}
 \label{Expression_phi_-} \phi_-(s)&=&\frac{(1-\rho)\delta^2(\omega,\lambda)}{\omega}\frac{1}{\delta(\omega,\lambda)-s}, ~~~  \rm {Re}\ s\leq0,\\
 \label{Expression_phi_+} \phi_+(s)&=&\frac{(1-\rho)\delta(\omega,\lambda)}{\omega} \left[\frac{(\delta(\omega,\lambda)-\omega)s-\lambda\delta(\omega,\lambda)(1-\beta(s))}{(s-\lambda(1-\beta(s)))(s-\delta(\omega,\lambda))}\right],
~~~ \rm {Re}\ s\geq0.
 \end{eqnarray}\label{final_exp_phi_+}
One imediately sees that the density $v_-(x)$ is exponential:
 \begin{equation}{\label{Expression_v_-}}
 v_-(x)=\frac{(1-\rho)\delta^2(\omega,\lambda)}{\omega} e^{-\delta(\omega,\lambda)x} \ \text{for}\ x >0.
 \end{equation}
 One can rewrite the expression for $\phi_+(s)$ in the following form:
 \begin{equation}
 \phi_+(s)=\left[(1-\rho)\lambda\frac{(1-\beta(s))}{s-\lambda(1-\beta(s))}\right]\left[\frac{\delta(\omega,\lambda)}{\omega}
 \frac{(\frac{\delta(\omega,\lambda)-\omega}{\lambda})\frac{s}{1-\beta(s)}-\delta(\omega,\lambda)}{s-\delta(\omega,\lambda)}\right].
\label{23}
 \end{equation}
{\bf Remark}.
The previous formula expresses the function $\phi_+(s)$ as a product of the Laplace transform of the density of the M/G/1 workload,
viz., $\frac{(1-\rho) \lambda (1-\beta(s))}{s-\lambda(1-\beta(s))}$ and a second factor.
That fact has led us to the observation that
the queue behaves like an M/G/1 queue with different first service time in a busy period (cf.\ Welch \cite{Welch}
or pages 392-394 and 401 of Wolff \cite{Wolff}).
Let us explain this in some detail.

When restricting ourselves to the time intervals with a positive workload in the queue,
$v_+(x)$ behaves like the workload density in an $M/G/1$ queue with Poisson($\lambda$) arrival process and with i.i.d.\ service times
$B_1,B_2,\dots$ with distribution $B(\cdot)$ and LST $\beta(\cdot)$, but with the first service time $\hat{B}$
of each busy period having a {\em different} distribution $\hat{B}(\cdot)$ with LST $\hat{\beta}(\cdot)$.
Indeed, $\hat{B}$ is distributed like the overshoot
above $0$ of a service time $B$, when starting from some negative value $-V_- = -x$
which, by PASTA, has steady-state density $v_-(x)$.
In fact, it is easy to determine $\hat{\beta}(\cdot)$, since $v_-(x)$ is exponentially distributed with parameter $\delta(\omega,\lambda)$, as seen in
(\ref{Expression_v_-}).
For simplicity of notation, we write $\delta := \delta(\omega,\lambda)$. Then
\begin{eqnarray}
\hat{\beta}(s) &=& \int_{0+}^{\infty} {\rm e}^{-sx} {\rm d}_x \mathbb P (B-V_- <x|B-V_- > 0)
\nonumber
\\
&=& \frac{
\int_{0+}^{\infty} {\rm e}^{-sx} {\rm d}_x \mathbb P (B-V_- <x)}{ \mathbb P (B-V_- > 0)}
\nonumber
\\
&=& \frac{
\int_{x=0+}^{\infty} {\rm e}^{-sx} \int_{y=0}^{\infty} \delta {\rm e}^{-\delta y}  {\rm d}_x \mathbb P (B<x+y) {\rm d}y}{1-\beta{\delta}}
\nonumber
\\
&=&
\frac{\delta}{s-\delta} \frac{\beta(\delta)-\beta(s)}{1-\beta(\delta)} ,
\end{eqnarray}
and
\begin{equation}
1-\hat{\beta}(s) = \frac{s}{s-\delta} - \frac{\delta}{s- \delta} \frac{1-\beta(s)}{1-\beta(\delta)} .
\label{1minushat}
\end{equation}
From \cite{Welch} or p.\ 401 of \cite{Wolff} it is seen that the LST of the steady-state workload (and waiting time) distribution
in this queueing model with exceptional first service
is given by
\begin{equation}
\mathbb E [{\rm e}^{-sW}] = \pi_0 \lambda \frac{\beta(s) - \hat{\beta}(s)(1-s/\lambda)}{s-\lambda(1-\beta(s))},
\end{equation}
with $\pi_0$ the probability of an empty system,
and hence
\begin{equation}
\mathbb E [{\rm e}^{-sW}|W>0] = \frac{\mathbb E [{\rm e}^{-sW}] - \pi_0}{1-\pi_0} = \frac{\pi_0 \lambda}{1-\pi_0} \frac{1 - \hat{\beta}(s)}{s-\lambda(1-\beta(s))}.
\label{W|W>0}
\end{equation}
A balance argument, or the observation that the expressions in (\ref{W|W>0}) should equal $1$ for $s=0$,
and that the Laplace-Stieltjes transforms
of the residual ordinary service time
$\frac{1-\beta(s)}{s \mathbb E B}$
and of the residual special service time
$\frac{1-\hat{\beta}(s)}{s \mathbb E \hat{B}}$ equal $1$ for $s=0$, readily yields that
\begin{equation}
\pi_0 = \frac{1-\rho}{1 -\rho + \lambda \mathbb E \hat{B}} .
\end{equation}
It readily follows from (\ref{1minushat}) that here
\begin{equation}
\mathbb E \hat{B} = \frac{\mathbb E B}{1 - \beta(\delta)} - \frac{1}{\delta} .
\end{equation}
Now compare (\ref{W|W>0}) and (\ref{Expression_phi_+}).
We claim they agree up to a multiplicative constant, which is $\int_0^{\infty} v_+(x) {\rm d}x$.
Indeed, one can rewrite the term between square brackets in (\ref{Expression_phi_+}) as
(replace $\delta-\omega$ by $\lambda (1-\beta(\delta))$, using the fact that
$1-\beta(\delta) = \frac{\delta - \omega}{\lambda}$,
which follows from the definition of $\delta$ as the zero of $s-\lambda(1-\beta(s)) = \omega$),
\begin{equation}
(1-\beta(\delta)) [\frac{1}{s-\lambda(1-\beta(s))} (\frac{s}{s-\delta}
- \frac{\delta}{1-\beta(\delta)} \frac{1-\beta(s)}{s-\delta})] .
\label{tussen}
\end{equation}
Now use relation (\ref{1minushat})
to see that the factor between square brackets in (\ref{tussen}) equals
the factor
$\frac{1 - \hat{\beta}(s)}{s-\lambda(1-\beta(s))}$ in (\ref{W|W>0}).
\\

\noindent
{\bf Remark}.
It may at first sight seem surprising that $v_-(x)$ is an exponential density.
The fact that $v_-(x)$ is exponential may be explained as follows.
Consider the inventory process, so look at Figure $2$ upside down (``process $1$").
Next, consider this figure (now reproduced as the top figure in Figure $3$) by looking from right to left (``process $2$").
Subsequently, replace the line segments
that go up at an angle of 45 degrees by upwards jumps equal to the increase along the line segment; and replace the jumps downward by line segments that go down at an angle of 45 degrees, by an amount equal to the jump (``process $3$"; see the bottom figure in Figure $3$).
One now has the representation of the workload process in a busy period of a $G/M/1$ queue.
Indeed, the jumps upward (service times) are exp($\lambda$) distributed,
and the intervals between jumps have distribution $B(\cdot)$.
Notice in particular that the waiting times in the $G/M/1$ queue are identical to the heights after jumps
in process $2$. By PASTA,
these heights have the same distribution as the steady-state workload distribution in process $2$, and hence also in process $1$.
Finally use the fact that the waiting time in the $G/M/1$ queue is exponentially distributed.
\begin{figure}[!h] \label{Fig3}
\begin{center}
\includegraphics[width=12cm]{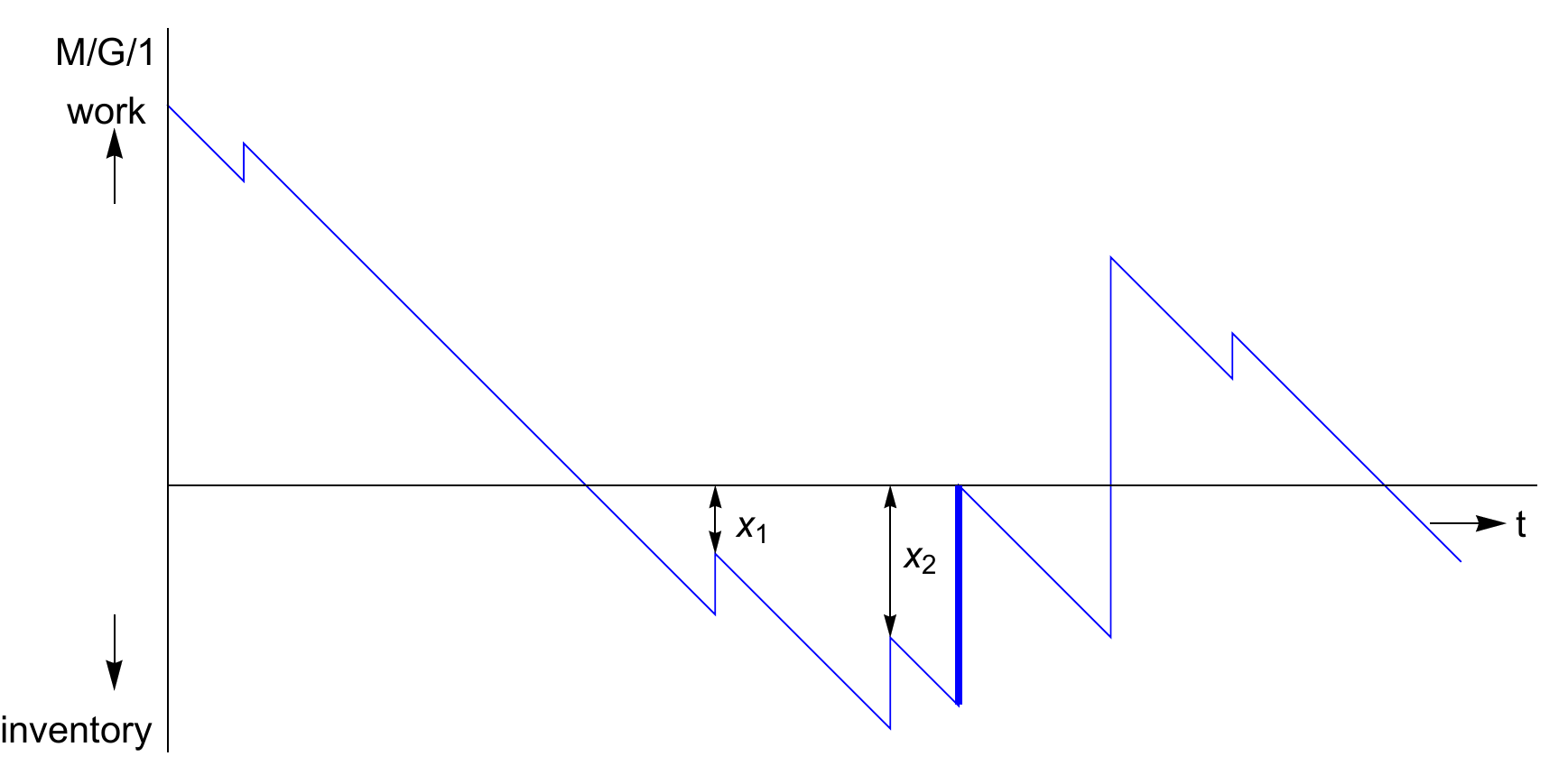}\\[3mm]
\caption{The workload and inventory process}
\includegraphics[width=12cm]{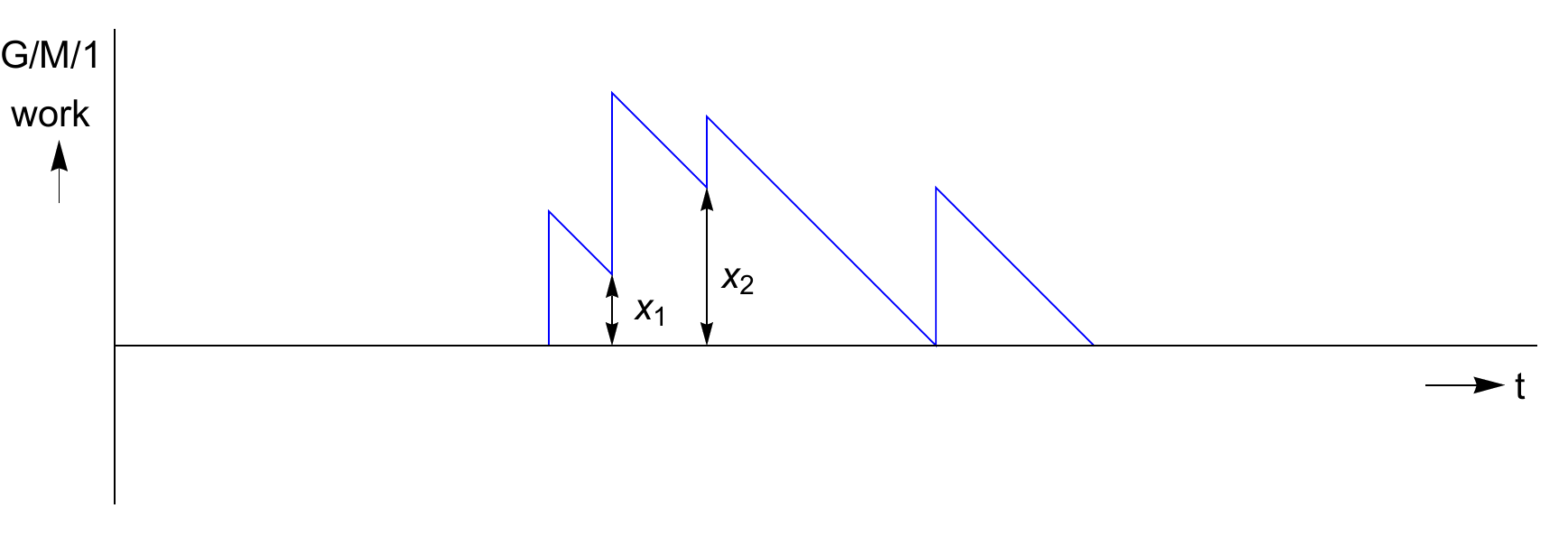}
 \caption{$G/M/1$ busy period (``process $3$")}
\end{center}
\end{figure}
\begin{example}{Exponential service requirements in the queueing/inventory model}\\
 We will keep the same notations as previously.
 In this case, one has $\mathbb P(B>x)=e^{-\mu x}$ with $\mu >0$ and
 for $\rm {Re}\ s\geq0$, $\beta(s)=\frac{\mu}{\mu+s}$;
$\mathbb E(B)=\frac{1}{\mu}$, $\rho=\frac{\lambda}{\mu}$.
 The functions $f_{\omega,\lambda}$ and $f_{0,\lambda} $ in this case are given by
 \begin{equation}{\label{expression_f(omega,lambda)}}
 f_{\omega,\lambda}(s)=\frac{s^2+(\mu-\lambda-\omega)s-\omega\mu}{\mu+s}, \quad\text{and}\ ~~~ f_{0,\lambda}(s)=\frac{s(s+\mu-\lambda)}{s+\mu}.
 \end{equation}
 The function $f_{\omega,\lambda}$ has two zeros:
 \begin{eqnarray*}
  \delta(\omega, \lambda)&=&\frac{\sqrt{(\mu-\lambda-\omega)^2+4\omega\mu}-(\mu-\lambda-\omega)}{2}>0,\\
  \text{and}\ \eta(\omega, \lambda)&=&\frac{-\sqrt{(\mu-\lambda-\omega)^2+4\omega\mu}-(\mu-\lambda-\omega)}{2}<0.
 \end{eqnarray*}
 The functions $g_{\omega,\lambda}$ and $g_{0,\lambda}$ are then given by
 \begin{equation}{\label{expression_g(omega,lambda)}}
 g_{\omega,\lambda}(s)=\frac{s-\eta(\omega,\lambda)}{s+\mu}, \quad\text{and}\ ~~~ g_{0,\lambda}(s)=\frac{s+\mu-\lambda}{s+\mu}.
 \end{equation}
$A$ is given by Equation ({\ref{Calcul_A}}).
The function $v_-(x)$ for $x>0$ is given by ({\ref{Expression_v_-}}), so:
\begin{eqnarray}{\label{Expression_v_-_exp}}
v_-(x)&=&\frac{\mu-\lambda}{2\omega}\left[(\mu-\lambda-\omega)^2-2\omega\mu-(\mu-\lambda-\omega)\sqrt{(\mu-\lambda-\omega)^2+4\omega\mu}\right]\\
&\ &\times \exp\left[-\frac{\sqrt{(\mu-\lambda-\omega)^2+4\omega\mu}-(\mu-\lambda-\omega)}{2}x\right]. \nonumber
\end{eqnarray}
Using Relations ({\ref{factors_identification_ph_+_1}}) and
({\ref{expression_g(omega,lambda)}}), and after some calculations, one gets for $\rm {Re}\ s\geq0$:
\begin{equation}
\phi_+(s)=\frac{(\mu-\lambda)(\mu-\eta(\omega,\lambda))\delta^2(\omega,\lambda)}{\mu^2\omega}\frac{1}{s+\mu-\lambda}.
\end{equation}
Consequently, one can deduce for $x>0$:
\begin{equation}
v_+(x)=\frac{ (\mu-\lambda)(\mu-\eta(\omega,\lambda))\delta^2(\omega,\lambda)}{\mu^2\omega}e^{-(\mu-\lambda)x}.
\end{equation}
The following figures represent respectively the steady state inventory density $v_-$ and the steady state workload density
$v_+$, in the exponential service requirements case for the particular values $\mu=2,\ \lambda=1$ and $\omega=2$.
\begin{figure}[!h]
\begin{center}
\includegraphics[width=0.7\textwidth]{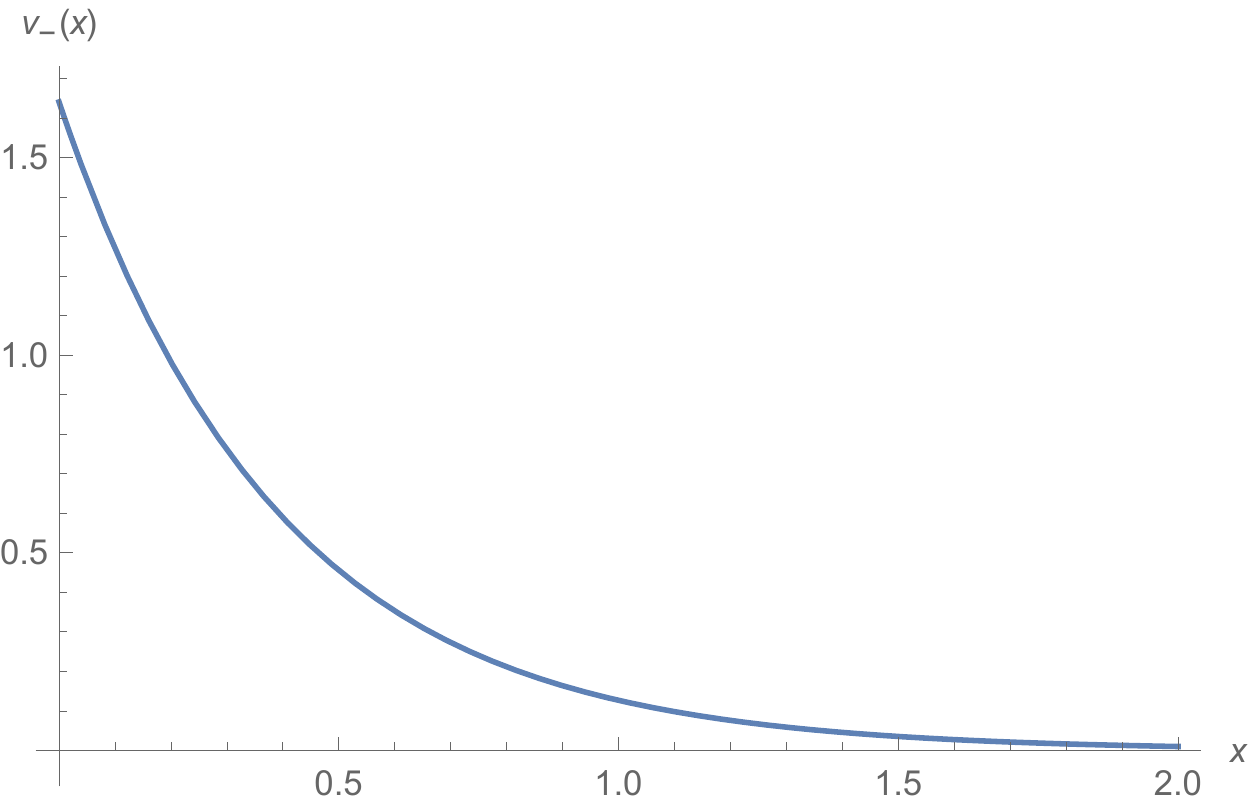}
\caption{Steady state inventory density ($\mu=2,\ \lambda=1$ and $\omega=2$).}
\end{center}
\end{figure}
\begin{figure}[!h]
\begin{center}
\includegraphics[width=0.7\textwidth]{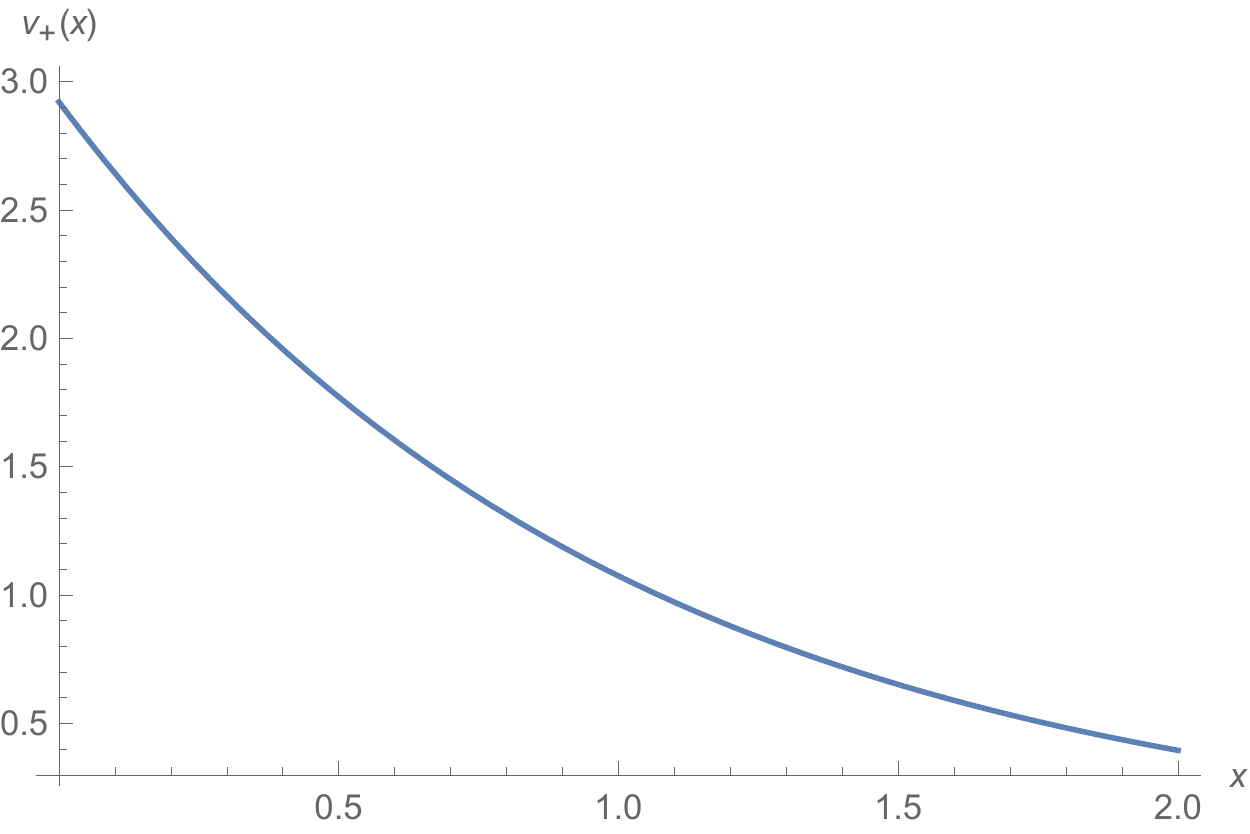}
\caption{Steady state workload density ($\mu=2,\ \lambda=1$ and $\omega=2$).}
\end{center}
\end{figure}
\end{example}
\subsection{Insurance model}\ \\
We now turn to the insurance risk model with bankruptcy.
Equation (\ref{Wiener-Hopf_insurance_1}) becomes, for ${\rm Re}~s =0$:
\begin{equation} \label{Wiener-Hopf_insurance_2}
\left(cs - \lambda(1- \beta(s))\right)\Psi_+(s)+\frac{\lambda(1-\beta(s))-cs}{s} =\left(cs - \lambda (1- \beta(s))-\omega\right)\Psi_-(-s) .
\end{equation}
In particular, when $s=0$ in (\ref{Wiener-Hopf_insurance_2}), one has
$\displaystyle \lambda \lim_{s\rightarrow 0} \frac{1-\beta(s)}{s}-c=-\omega\Psi_-(0)$, which implies
\begin{equation}\label{Identification_Psi-(0)}
\Psi_-(0)=\frac{c-\lambda \mathbb E(Y)}{\omega}.
\end{equation}
We follow the same procedure as in the queueing/inventory model: We reformulate Equation (\ref{Wiener-Hopf_insurance_2}) into a Wiener-Hopf
problem. Set
$f_{\omega,\lambda,c}:\quad s\mapsto cs - \lambda (1- \beta(s))-\omega,\  {\rm{Re}}\ s\geq0$ and $f_{0,\lambda,c}:\quad
s\mapsto cs - \lambda (1- \beta(s)),\  {\rm{Re}}\ s\geq0$.
According to \cite{Cohen} p.\ $548$, the constant $\omega$ being positive, the function $f_{\omega,\lambda,c}$ has one zero
$s=\delta(\omega,\lambda,c)$ and this zero is simple satisfying ${\rm{Re}}\ \delta(\omega,\lambda,c)>0$. Also, the function $f_{0,\lambda,c}$ has $s=0$ as its only zero and this zero is simple. In particular the functions
$g_{\omega,\lambda,c}:\quad s \mapsto \frac {f_{\omega,\lambda,c}(s)}{s-\delta(\omega,\lambda,c)}$
and $g_{0,\lambda,c}: s \mapsto \frac {f_{0,\lambda,c}(s)}{s}$ are analytic for ${\rm{Re}}\ s>0$, continuous for ${\rm{Re}}\ s\geq0$ and take non-zero values on ${\rm{Re}}\ s\geq0$.
Dividing by $cs -\lambda (1- \beta(s))-\omega$ and multiplying by $s-\delta(\omega,\lambda,c)$
in (\ref{Wiener-Hopf_insurance_2}), one gets for ${\rm{Re}}\ s=0$:
\begin{eqnarray} \label{Wiener-Hopf_insurance_3}
&\frac{cs - \lambda(1- \beta(s))}{cs -\lambda (1-\beta(s))-\omega}\left(s-\delta(\omega,\lambda,c)\right)\Psi_+(s)
+ \frac{\lambda(1-\beta(s))-cs}{s} \frac{s-\delta(\omega,\lambda,c)}{cs -\lambda (1-\beta(s))-\omega}\\
&=\Psi_-(-s)\left(s-\delta(\omega,\lambda,c)\right) \nonumber.
\end{eqnarray}
Clearly, the lefthand side of (\ref{Wiener-Hopf_insurance_3}) is analytic for ${\rm{Re}}\ s>0$ and continuous for ${\rm{Re}}\ s\geq 0$;
on the other hand, the righthand side is analytic for ${\rm{Re}}\ s<0$ and continuous for ${\rm{Re}}\ s\leq 0$. In addition, both sides coincide
for ${\rm{Re}}\ s=0$. By Liouville's theorem, there exist $n\geq 0$ and a polynomial $R_n$ of degree $n$ such that
\begin{eqnarray}
&\frac{cs -\lambda(1- \beta(s))}{cs -\lambda (1-\beta(s))-\omega}\left(s-\delta(\omega,\lambda,c)\right)\Psi_+(s)
\label{Identification_Psi+(s)}\\
\vspace{0.07cm}
&+ \frac{\lambda(1-\beta(s))-cs}{s}\frac{s-\delta(\omega,\lambda,c)}{cs-\lambda (1-\beta(s))-\omega}=R_n(s)\ \text{for}\ {\rm{Re}}\ s\geq 0, \nonumber\\
&\Psi_-(-s)\left(s-\delta(\omega,\lambda,c)\right)=R_n(s)\ \text{for}\ {\rm{Re}}\ s\leq 0 \label{Identification_Psi-(s)}.
\end{eqnarray}
Since $\displaystyle \lim_{s \rightarrow -\infty} \Psi_-(-s)=0$ and using Equation (\ref{Identification_Psi-(s)}),
we can deduce that $n$ must be $0$; say $R_n(s)=Z, \ Z\in \mathbb C$.
Consequently, one has
\begin{eqnarray}
\Psi_+(s)&=&Z\frac{cs-\lambda (1-\beta(s))-\omega}{s-\delta(\omega,\lambda,c)}\frac{1}{cs-\lambda (1-\beta(s))}
+\frac{1}{s},
~~~ \text{for}\ {\rm{Re}}\ s\geq 0,
\label{Identification_Psi+(s)_1}\\
\Psi_-(-s)&=&\frac{Z}{s-\delta(\omega,\lambda,c)}, ~~~ \text{for}\ {\rm{Re}}\ s\leq 0 \label{Identification_Psi-(s)_1}.
\end{eqnarray}
In particular, Equation (\ref{Identification_Psi-(s)_1}) implies
\begin{equation}\label{Identification_Psi-(0)}
\Psi_-(0)=\frac{-Z}{\delta(\omega,\lambda,c)}.
\end{equation}
Let us now identify the constant $Z$.\\
Set $\rho=\frac{\lambda \mathbb E(Y)}{c}.$
Plugging $s=0$ in Equation (\ref{Wiener-Hopf_insurance_2}) and combining it with Relation (\ref{Identification_Psi-(0)}), one gets
 \begin{equation}\label{Identification_A}
 Z=-\frac{c\delta(\omega,\lambda,c)}{\omega}\left(1-\rho\right) ;
\end{equation}
thanks to Equation (\ref{Identification_Psi-(s)_1}), the latter relation completely determines the function $\Psi_-(-s)$ for ${\rm{Re}}\ s\leq 0$.
In this case, one can immediately deduce
$u_-(x)$, which is the survival probability when starting at a negative surplus $-x$, and so also the ruin probability $\tilde{u}_-(x)=1-u_-(x)$. One has for $x\geq 0$,
\begin{equation}\label{Identification_u_-}
\tilde{u}_-(x) = 1 - u_-(x)= 1 - \frac{c\delta(\omega,\lambda,c)(1-\rho)}{\omega}e^{-\delta(\omega,\lambda,c)x} .
\end{equation}
Finally, since the constant $Z$ is known, we can identify the function $\Psi_+(s)$ for ${\rm{Re}}\ s\geq 0$. Rewriting Equation (\ref{Identification_Psi+(s)_1}), one has
\begin{equation}\label{Identification_Psi+(s)_2}
\Psi_+(s)=\frac{1}{g_{0,\lambda,c}(s)}\frac{Zg_{\omega,\lambda,c}(s)-\lambda \frac{1-\beta(s)}{s}+c}{s}\  \text{for}\ {\rm{Re}}\ s\geq 0.
\end{equation}
Set $h_{\omega,\lambda,c}:\quad s\mapsto Zg_{\omega,\lambda,c}(s)-\lambda \frac{1-\beta(s)}{s}, \ \text{for} \ {\rm{Re}}\ s\geq 0$,
so $h_{\omega,\lambda,c}(0)= \frac{Z\omega}{\delta(\omega,\lambda,c)}-\lambda \mathbb E(Y)$. Equation
(\ref{Identification_A}) implies that $h_{\omega,\lambda,c}(0)=-c$. Then one has
\begin{equation}\label{Identification_Psi+(s)_3}
\Psi_+(s)=\frac{1}{g_{0,\lambda,c}(s)}\frac{h_{\omega,\lambda,c}(s)-h_{\omega,\lambda,c}(0)}{s} \ ;
\end{equation}
 the function $s\mapsto\frac{h_{\omega,\lambda,c}(s)-h_{\omega,\lambda,c}(0)}{s}$ being analytic for ${\rm{Re}}\ s>0$, continuous for ${\rm{Re}}\ s\geq0$
 with $h_{\omega,\lambda,c}(0)=-c$.
\\
The Laplace transform of the survival probability when starting at a positive surplus $x$, given by the function $1-u_+(x)$,
equals:
\begin{equation}
\frac{1}{s}-\Psi_+(s) = Z \frac{1}{cs-\lambda(1-\beta(s))} \frac{cs -\lambda(1- \beta(s)) -\omega}{\delta(\omega,\lambda,c) - s} .
\end{equation}

\begin{example}{Exponential claim sizes in the insurance model}\\
We will keep the same notations as previously. Fix $ \nu>0$.
In the exponential claim sizes case, we assume claim size density
$\nu e^{-\nu y}$, and hence $$\beta(s)=\frac{\nu}{s+\nu}\ \text{for}\ {\rm{Re}}\ s\geq0.$$
In particular $\rho=\frac{\lambda}{\nu c}$.
In this case, one obtains
\begin{equation}\label{f_Exp_case}
f_{\omega,\lambda,c}=\frac{cs^2+(c\nu-\omega-\lambda)s-\omega\nu}{s+\nu} \quad \text{for}\ {\rm{Re}}\ s\geq0.
\end{equation}
This function has two zeros:
$$ \eta(\omega,\lambda,c)=\frac{-\sqrt{(c\nu-\omega-\lambda)^2+4\omega\nu c}-(c\nu-\omega-\lambda)}{2c} <0,$$
$$\text{and}\ \delta(\omega,\lambda,c)=\frac{\sqrt{(c\nu-\omega-\lambda)^2+4\omega\nu c}-(c\nu-\omega-\lambda)}{2c} >0.$$ Therefore,
\begin{equation}\label{g_Exp_case}
g_{\omega,\lambda,c}=\frac{f_{\omega,\lambda,c}(s)}{s-\delta(\omega,\lambda,c)}=c\frac{s-\eta(\omega,\lambda,c)}{s+\nu} \quad \text{for}\ {\rm{Re}}\ s\geq0.
\end{equation}
Applying Equation (\ref{Identification_A}), one gets
\begin{equation}\label{identification_A_Exp}
Z=-(\nu-\frac{\lambda}{c})\frac{c}{\nu\omega}\delta(\omega,\lambda,c).
\end{equation}
Using the relation between the zeros $\delta(\omega,\lambda,c)$ and $\eta(\omega,\lambda,c)$,
$$\delta(\omega,\lambda,c)\eta(\omega,\lambda,c)=-\frac{\omega\nu}{c},$$ one obtains
\begin{equation}\label{identification_A_Exp_1}
Z=\frac{\nu-\frac{\lambda}{c}}{\eta(\omega,\lambda,c)}.
\end{equation}
Applying Equation (\ref{Identification_u_-}), one gets for $x\geq 0$,
\begin{equation}\label{Identification_u_-_Exp}
u_-(x)=\frac{\nu-\frac{\lambda}{c}}{-\eta(\omega,\lambda,c)}e^{-\delta(\omega,\lambda,c)x}\ \text{and}\
\tilde u_-(x)=1-\frac{\nu-\frac{\lambda}{c}}{-\eta(\omega,\lambda,c)}e^{-\delta(\omega,\lambda,c)x}.
\end{equation}
These results agree with results of Albrecher and Lautscham in \cite{AlbL}, Section 2.1.1.
To explicitly determine the function $ \Psi_+$ as given in
(\ref{Identification_Psi+(s)_3}),
and so $u_+$, one has to make
the function $h_{\omega,\lambda,c}$ explicit in this case. After some calculations, one gets for ${\rm{Re}}\ s\geq0$,
 \begin{equation}\label{Calcul_h_Exp}
 h_{\omega,\lambda,c}(s)=\frac{(\nu c-\lambda)s-\nu c \eta(\omega,\lambda,c)}{\eta(\omega,\lambda,c)(s+\nu)}.
 \end{equation}
 One now can use Equation (\ref{Identification_Psi+(s)_3}) and deduce for ${\rm{Re}}\ s\geq0$
 \begin{equation}\label{Identification_Psi+(s)_Exp}
 \Psi_+(s)=\left(1-\frac{\nu-\frac{\lambda}{c}}{-\eta(\omega,\lambda,c)}\right)\frac{1}{s+\nu-\frac{\lambda}{c}}.
 \end{equation}
 Notice here that the condition $\rho<1$ is equivalent to $\nu-\frac{\lambda}{c}>0$.
 Finally, one obtains for $x\geq 0$:
 \begin{equation}\label{Calcul_u_+_Exp}
 u_+(x)=\left(1-\frac{\nu-\frac{\lambda}{c}}{-\eta(\omega,\lambda,c)}\right)e^{-\left(\nu-\frac{\lambda}{c}\right)x}.
 \end{equation}
Note that using equations  (\ref{Identification_u_-_Exp}) and (\ref{Calcul_u_+_Exp}), one can check immediately that $u_+(0)+u_-(0)=1$.
Furthermore, note that
(\ref{Identification_u_-_Exp}) and (\ref{Calcul_u_+_Exp})
coincide with Formula (18) in \cite{AlbL}.\\
The following figures represent respectively the bankruptcy probability starting from a negative surplus against the initial surplus $-C_0$
and the bankruptcy probability starting from a positive surplus against the initial surplus $C_0$, in the exponential claim sizes case for the particular values $\nu=2,\ \lambda=1,\ \omega=2$ and $c=1$.
\setcounter{topnumber}{3}
\begin{figure}[!h]
\begin{center}
\includegraphics[width=0.7\textwidth]{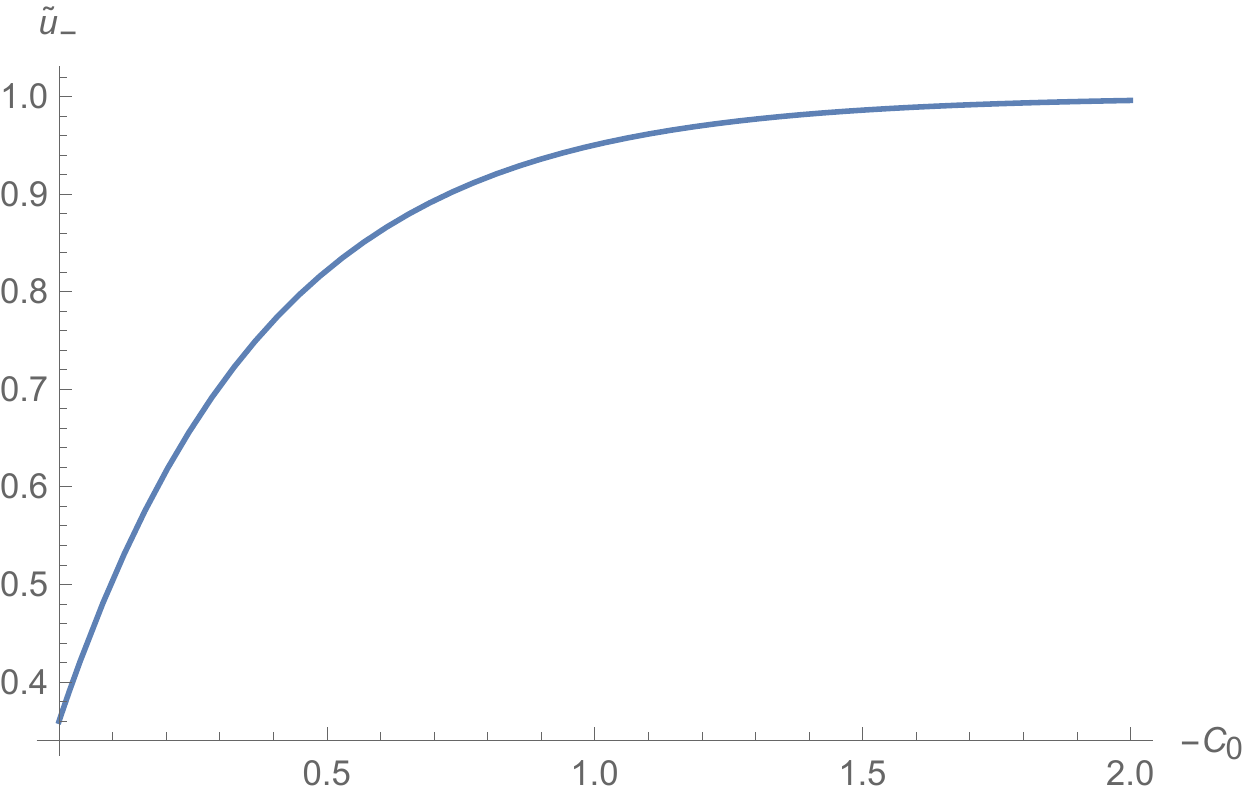}
\caption{Bankruptcy probability starting from a negative surplus against the initial surplus $C_0$ ($\nu=2,\ \lambda=1,\ \omega=2$ and $c=1$).}
\end{center}
\end{figure}\ \\
\begin{figure}[!h]
\begin{center}
\includegraphics[width=0.7\textwidth]{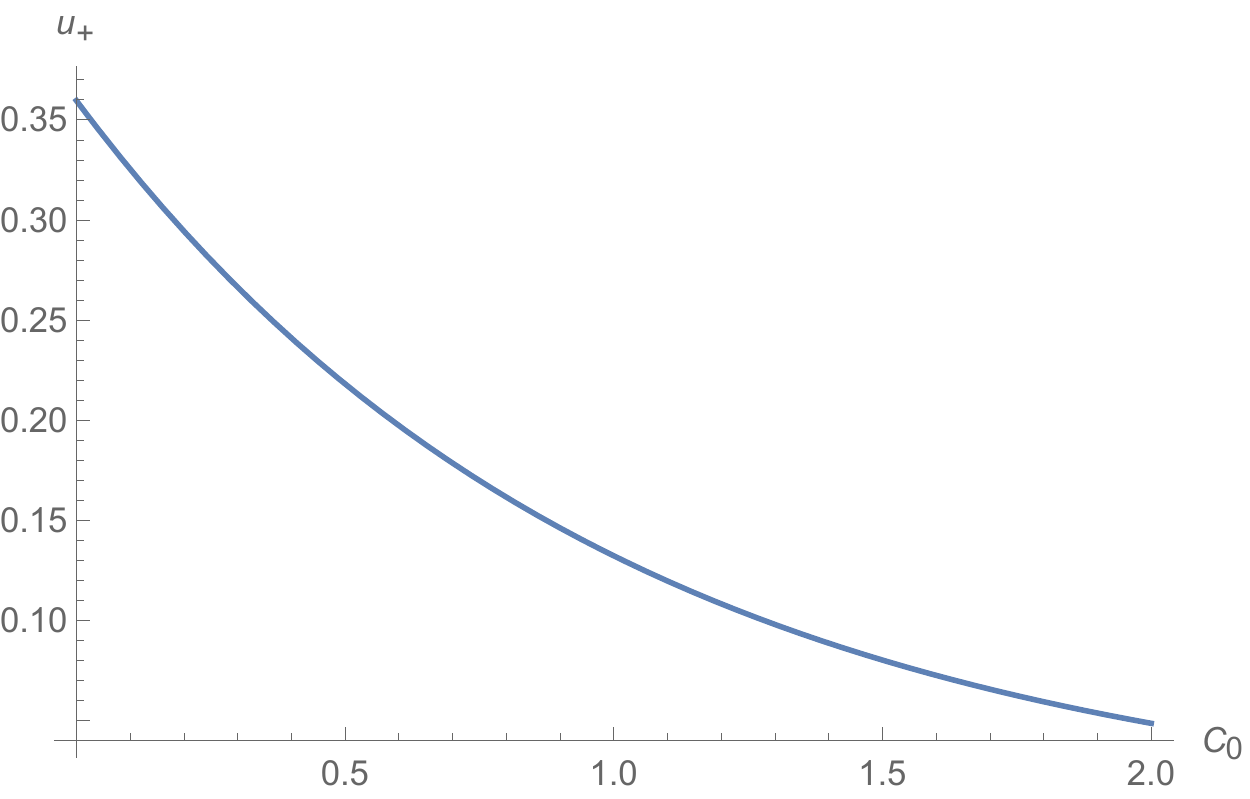}
\caption{Bankruptcy probability starting from a positive surplus against the initial surplus $C_0$ ($\nu=2,\ \lambda=1,\ \omega=2$ and $c=1$).}
\end{center}
\end{figure}
\end{example}\ \\
\vspace{0.15cm}\ \\

\begin{remark}
The queueing/inventory and insurance risk models that we treat in this paper are clearly closely related,
although they are not dual in the sense discussed in, e.g., Section III.2 of \cite{AA}.
The results that we obtain for the densities and their Laplace transforms are indeed very similar; cf.\ Equations
(\ref{Expression_phi_+}) and
(\ref{Identification_Psi+(s)_2}),
and (\ref{Expression_v_-}) and (\ref{Identification_u_-}).
It would be interesting to construct a queueing/inventory model that is completely
dual to the insurance risk model.
\end{remark}
\begin{remark}
Albrecher and Ivanovs \cite{AI} have recently studied exit problems for L\'evy processes where the first passage time over a threshold is detected either
immediately ('ruin') or at Poisson observation epochs ('bankruptcy').
They relate the two exit problems via a nice identity.
In the case of the Cram\'er-Lundberg insurance risk model, their identity reads as follows:
\begin{equation}
\hat{s}(x) = \mathbb E[s(x+U)],
\label{AIhat}
\end{equation}
where $s(x)$ is the survival probability in the case of the Cram\'er-Lundberg model with initial capital $x$
and $\hat{s}(x)$ is the survival probability in the corresponding model
with Poisson($\omega$) observations (notice that $u_+(x) = 1 - \hat{s}(x)$);
finally, $U$ is an exp($\Phi$) distributed random variable, where $\Phi$ is the inverse of the Laplace exponent
of the spectrally negative L\'evy process corresponding to the Cram\'er-Lundberg model. In other words, $\Phi$
is the zero of $cs - \lambda(1-\beta(s)) = \omega$. We conclude that $\Phi = \delta(\omega,\lambda,c)$.
\\
It follows from (\ref{AIhat}) that
\begin{equation}
\hat{s}(x) = \int_{t=0}^{\infty} \delta(\omega,\lambda,c) {\rm e}^{- \delta(\omega,\lambda,c)t} s(x+t) {\rm d}t .
\label{AIhat2}
\end{equation}
When $s(\cdot)$ is explicitly known, then one can determine $\hat{s}(x)$ explicitly
using (\ref{AIhat2}).
In particular, in the case of exp($\nu$) claim sizes, one has (cf.\ \cite{AA})
\begin{equation}
s(x) = 1 - \frac{\lambda}{\nu c} {\rm e}^{-(\nu - \frac{\lambda}{c})x}, ~~~~ x > 0,
\end{equation}
and hence
\begin{equation}
\hat{s}(x) = 1 - \frac{\delta(\omega,\lambda,c)}{\delta(\omega,\lambda,c)+\nu c - \lambda}
\frac{\lambda}{\nu c} {\rm e}^{-(\nu - \frac{\lambda}{c})x}, ~~~ x > 0.
\end{equation}
Using the definition of $\delta(\omega,\lambda,c)$, it is readily verified that this formula
is indeed in agreement with the expression for $u_+(x) = 1 - \hat{s}(x)$ in
(\ref{Calcul_u_+_Exp}).
\end{remark}

\section{Analysis for $\omega(\cdot)$ linear}
\label{linearomega}
In this section we focus on the queueing/inventory model.
We assume in this section that the function $\omega$ introduced in Section~\ref{Qmodel} is linear, i.e., there exists a constant $a>0$ such that
for all $x\geq 0$, one has
\begin{equation}\label{omega_linear}
\omega(x)=ax.
\end{equation}
Equation ({\ref{eq_Wiener_Hopf_0}}) becomes
\begin{equation}{\label{eq_Wiener_Hopf_lin_1}}
\left[s-\lambda
(1-\beta(s))\right]\phi_+(s)=-\left[s-\lambda(1-\beta(s))\right]\phi_-(-s)-a\int_0^{+\infty}y(1-e^{sy})v_-(y)dy
 , \quad \text{for}
\quad \rm {Re}\ s=0.
\end{equation}
Set $\mathbb E I=\int_0^{+\infty}yv_-(y)dy$.
After integrating by parts, one gets
\begin{equation}{\label{eq_Wiener_Hopf_lin_2}}
\left[s-\lambda
(1-\beta(s))\right]\phi_+(s)=-\left[s-\lambda(1-\beta(s))\right]\phi_-(-s)+a\frac{d}{ds}[\phi_-(-s)]-a \mathbb E I
 , \quad \text{for}
\quad \rm {Re}\ s=0.
\end{equation}
We will discuss here the case where the function $\beta$ is rational: suppose that there exist $m \in \mathbb N$ and polynomials $N$ and $D$ in $\mathbb C[x]$ such that ${\rm deg}D=m$, ${\rm deg}N\leq m-1$, $(N\bigvee D)=1$, i.e they do not have a common factor and $\beta(s)=\frac{N(s)}{D(s)}\ \text{for}\ \rm {Re}\ s=0$. In this configuration, necessarily, the polynomial $D$ has no zeros for
$\rm {Re}\  s \geq 0$ and $m$ zeros for $\rm {Re}\  s \leq 0$ (counted with multiplicities). Denote them by $-\mu_1,-\mu_2,\cdots,-\mu_m$
with $\rm {Re}\  \mu_j \geq 0$.
Set $N(s)=\displaystyle \sum_{k=0}^{m-1}n_ks^k$ and $D(s)=\displaystyle \sum_{k=0}^{m}d_ks^k$
for $s\in \mathbb C$, where $n_0,\cdots,n_{m-1}$ and $d_0,\cdots,d_m$ are complex numbers. Notice that $\beta(0)=1$ implies that $n_0=d_0$.
Multiplying Equation ({\ref{eq_Wiener_Hopf_lin_2}}) by $D(s)$, one obtains for $\rm {Re}\ s=0$:
\begin{equation}{\label{eq_Wiener_Hopf_lin_3}}
\left[(s-\lambda)D(s)
+\lambda N(s)\right]\phi_+(s)=-\left[(s-\lambda)D(s)
+\lambda N(s)\right]\phi_-(-s)+aD(s)\frac{d}{ds}[\phi_-(-s)]-a \mathbb E I D(s).
\end{equation}
Using the same techniques as previously, one can deduce that there exists a polynomial
$R_m$ such that
\begin{eqnarray}
\left[(s-\lambda)D(s)
+\lambda N(s)\right]\phi_+(s)&=R_m(s) \quad \text{for}
\quad \rm {Re}\ s\geq 0, \label{identification_phi_+_lin_1}\\
-\left[(s-\lambda)D(s)
+\lambda N(s)\right]\phi_-(-s)+aD(s)\frac{d}{ds}[\phi_-(-s)]-a \mathbb E I D(s)&=R_m(s)\quad \text{for}
\quad \rm {Re}\ s\leq0. \label{identification_phi_-_lin_1}
\end{eqnarray}
Since $\displaystyle \lim_{s\rightarrow +\infty}\phi_+(s)=0$ and using Equation (\ref{identification_phi_+_lin_1}), one can deduce that
${\rm deg}R_m \leq m$.
Equation (\ref{identification_phi_+_lin_1}) implies that $R_m(0)=0$. Hence set $R_m(s)=\displaystyle \sum_{k=1}^{m}r_ks^k$ for $s\in \mathbb C$, where $r_1,\cdots,r_m\in \mathbb C$.\\
It follows from (\ref{identification_phi_+_lin_1}) that
\begin{equation}
\phi_+(s) = \frac{R_m(s)}{(s-\lambda) D(s) + \lambda N(s)}, ~~~ {\rm Re}~s \geq 0.
\label{phi+s}
\end{equation}
Using Equation (\ref{phi+s}), one can also express $\phi_+(0)$ in $r_1$.
Indeed
\begin{equation}\label{Phi_+(0)}
\phi_+(0)=\lim_{s\rightarrow 0} \frac{R_m(s)}{(s-\lambda)D(s)
+\lambda N(s)}=\frac{r_1}{d_0+\lambda(n_1-d_1)}.
\end{equation}
Since $\phi_+(0)+\phi_-(0)=1$, one can also deduce
\begin{equation}\label{Phi_-(0)}
\phi_-(0)= 1-\frac{r_1}{d_0+\lambda(n_1-d_1)}.
\end{equation}
Notice that by plugging $-\mu_j$ for $j \in \{1,2,\cdots,m\}$ in Equation (\ref{identification_phi_-_lin_1}), one gets the $m$ relations
\begin{equation}\label{phi_-(mu_j)}
-\lambda N(-\mu_j)\phi_-(\mu_j)=R_m(-\mu_j), \quad 1\leq j\leq m.
\end{equation}
For $s\in \mathbb C$, introduce the following notations: $\check N(s):=N(-s)$, $\check D(s):=D(-s)$ and $\check R_m(s):=R_m(-s)$.
Now, putting $z=-s$ in Equation (\ref{identification_phi_-_lin_1}), the function $ \phi_-$ is a solution of the following
first order differential equation:
\begin{equation}{\label{Eq-diff_Phi_-_1}}
a\check D(z)\frac{d}{dz}[\phi_-(z)]+\left[-(z+\lambda)\check D(z)
+\lambda \check N(z))\right]\phi_-(z)+a \mathbb E I \check D(z)+\check R_m(z)=0\quad \text{for}
\quad \rm {Re}\ z\geq 0.
\end{equation}
One can write the previous equation in the following form:
\begin{equation}{\label{Eq-diff_Phi_-_1_without_singularity}}
\frac{d}{dz}[\phi_-(z)]=-\frac{\lambda \check N(z)\phi_-(z)+\check R_m(z)}{a\check D(z)}+\frac{z+\lambda}{a}\phi_-(z)
-\mathbb E I=0\quad \text{for} \quad \rm {Re}\ z\geq 0.
\end{equation}
Equation ({\ref{phi_-(mu_j)}})
implies that
the function $z\mapsto \frac{\lambda \check N(z) \phi_-(z)+\check R_m(z)}{\check D(z)}$
is analytic on $\rm {Re}\ z>0$ and continuous on $\rm {Re}\ z\geq0$.
Equation ({\ref{Eq-diff_Phi_-_1}}) is a first order algebraic differential equation on the complex half plane $\{\rm{Re}\ z> 0\}$.
It has m singularities $\mu_1,\mu_2, \dots , \mu_m$; Relation ({\ref{Eq-diff_Phi_-_1_without_singularity})} shows that
these singularities are removable.

In principle one can solve this first order algebraic differential equation.
However, the singularities give rise to several technical difficulties.
Below we consider the case $m=1$, i.e., $exp(\mu)$ distributed service requirements, yielding one singularity $z=\mu$.
In Subsection \ref{linear_exp} we sketch its analysis. We formally solve the differential equation (\ref{Eq-diff_Phi_-_1_without_singularity})
only for ${\rm Re}~z > \mu$, where $\mu$ is the singularity.
We also determine the two missing constants $EI$ and $r_1$.
In addition, we formally invert the Laplace transform $\phi_-(z)$ to find $v_-(x)$,
but this inversion is not considered in detail.

In the appendix we treat the solution of differential equation (\ref{Eq-diff_Phi_-_1_without_singularity}) in much more detail.
We consider both the case $0 \leq z \leq \mu$ (Subsection \ref{sub_sec_A_1}) and the case $z \geq \mu$ (Subsection \ref{sub_sec_A_2}),
and we determine the two missing constants $EI$ and $r_1$.
It turns out that the value of $\sigma := \frac{\mu \lambda}{a}$ plays a key role in the analysis.
One has to distinguish the case of non-integer $\sigma$ (consider $K-1 < \sigma < K$, for each $K=1,2,\dots$)
and the case of integer $\sigma = K$.
In the latter case, we also invert $\phi_-(z)$, finding an explicit expression for $v_-(x)$ in terms of Hermite polynomials.
We have included the appendix because it exposes both a series of technical difficulties
-- even for the case of just a single singularity $\mu$ --
in handling differential equation (\ref{Eq-diff_Phi_-_1_without_singularity})
as well as a mathematical approach to handle such difficulties.

In Section \ref{direct} we briefly outline a completely different approach to the problem of finding $v_-(x)$ and $v_+(x)$ for the case $\omega(x) = ax$ and
exp($\mu$) service requirements.
In that section we do not use Laplace transforms, but differentiate (\ref{queuing model<0}) twice to get a second order non-linear differential equation
in $v_-(x)$, and we differentiate (\ref{queuing model>0}) once to get a simple first order differential equation
in $v_+(x)$.
The latter equation is easily solved; the solution of the former differential equation
is expressed in hypergeometric functions.
For $\sigma = \frac{\mu \lambda}{a}$ integer, this hypergeometric function reduces to the Hermite polynomial found in
Subsection \ref{sub_sec_A_8}.

Finally, we should add that we do not see how the approach in Section \ref{direct} can be extended to the case of
an Erlang, hyperexponential or, more generally, phase-type service requirement distribution,
as such distributions will give rise to a higher-order non-linear differential equation for $v_-(x)$.
On the other hand, the approach of the appendix towards the differential equation (\ref{Eq-diff_Phi_-_1_without_singularity}) seems promising
for such service requirement distributions, even though they give rise to multiple
singularities $\mu_1,\dots,\mu_m$.
\subsection{The Exponential service requirements case}\label{linear_exp}
In this case, one has $\mathbb P(B>x)=e^{-\mu x}$ with $\mu >0$ and
 for $\rm {Re}\ s\geq0$, $\beta(s)=\frac{\mu}{\mu+s}$ ($N(s)=\mu$ and $D(s)=\mu+s$).
 Then $m=1$ and since $R_1(0)=0$, we obtain that $R_1(s)=r_1s$ where $r_1 \in \mathbb C$. Equation ({\ref{Eq-diff_Phi_-_1}}) becomes
\begin{equation}{\label{Eq-diff_Phi_-_1_exp}}
a (\mu-z)\frac{d}{dz}[\phi_-(z)]+z(z+\lambda-\mu)\phi_-(z)+a \mathbb E I(\mu-z)-r_1z=0\quad \text{for}
\quad \rm{Re}\ z\geq 0.
\end{equation}
Equation ({\ref{Eq-diff_Phi_-_1_exp}}) is a first order algebraic differential equation on the complex half plane $\{\rm{Re}\ z> 0\}$.
It has one singularity which is $\mu$. That makes the study of this differential equation more complicated; therefore, in this paper,
we dedicate an appendix to a further study concerning this differential equation.
Remembering that $\phi_-(\mu)=\frac{r_1}{\lambda}$ and writing Equation ({\ref{Eq-diff_Phi_-_1_exp}}) in the following form:
\begin{equation}{\label{Eq-diff_Phi_-_1_exp*}}
\frac{d}{dz}[\phi_-(z)]=\frac{\mu}{\lambda a} \frac{\phi_-(z)-\frac{r_1}{\lambda}}{z-\mu}+\frac{z+\lambda}{a}\phi_-(z)-\mathbb E I-\frac{r_1}{a},
\end{equation}
one deduces that this singularity is regular.
Solving Equation ({\ref{Eq-diff_Phi_-_1_exp}}) on $\{\rm{Re}\ z>\mu\}$, using the fact that $\displaystyle\lim_{z\rightarrow +\infty}\phi_-(z)=0$, one gets
\begin{equation}{\label{Sol_Phi_-exp_>mu}}
\phi_-(z)=z(\frac{z-\mu}{\mu})^{\frac{\lambda\mu}{a}}e^{\frac{z(z+2\lambda)}{2a}}
 \int_{1}^{+\infty}\left[\mathbb E I +\frac{r_1zt}{a(zt-\mu)}\right]\left(\frac{\mu}{zt-\mu}\right)^{\frac{\lambda\mu}{a}}
e^{\frac{-zt(zt+2\lambda)}{2a}} dt,
\end{equation}
(notice here that for every $\alpha\in \mathbb R$, we consider the principal value of the function $z \mapsto z^\alpha$).
For $\rm{Re}\ z>\mu$, we introduce $\displaystyle F_{\alpha,\beta}(z):=z(\frac{z-\mu}{\mu})^{\alpha}e^{\frac{z(z+2\lambda)}{2a}}\int_{1}^{+\infty}\left(\frac{\mu}{zt-\mu}\right)^{\beta}
e^{\frac{-zt(zt+2\lambda)}{2a}}dt$. The function $F_{\alpha,\beta}$ is analytic on
$\{\rm{Re}\ z>\mu\}$ and one can check by the l'Hopital rule that the function $F_{\alpha,\beta}$
can be analytically continued
in $\mu$ for $\alpha\geq \beta-1$. Now fixing $\alpha=\frac{\lambda \mu}{a}$, Equation (\ref{Sol_Phi_-exp_>mu}) can be written as follows:
\begin{equation}{\label{Sol_Phi_-exp_>mu_1}}
\phi_-(z)=(\mathbb E I+\frac{r_1}{a})F_{\alpha,\alpha}(z)+\frac{r_1}{a}F_{\alpha,\alpha+1}(z) \ \text{for}\ \rm{Re}\ z>\mu.
\end{equation}
We will denote by $\mathcal L^{-1}$ the inverse Laplace transform, and apply what is commonly known as Mellin inverse formula
or the Bromwich integral: let $\gamma$ be any real number such that $\gamma>\mu$, then one has
\begin{equation}
v_-(x)=\frac{1}{2i\pi}\int_{\gamma-i\infty}^{\gamma+i\infty}\phi_-(z)e^{zx}dz=\frac{1}{2\pi}\int_{-\infty}^{+\infty}\phi_-(\gamma+i\omega)e^{(\gamma+i\omega) x}d\omega\ \text{for}\ x>0.
\end{equation}
Since the constant $\gamma$ here is chosen to be larger than $\mu$, one can use the expression of $\phi_-$ given
in Equation ({\ref{Sol_Phi_-exp_>mu_1}}). Set $\displaystyle f_{\alpha,\beta}(x)=\frac{1}{2i\pi}\int_{\gamma-i\infty}^{\gamma+i\infty}F_{\alpha,\beta} (z)e^{zx}dz$. One gets
\begin{equation}{\label{v_-(x)_linear-exp}}
v_-(x)=(\mathbb E I +\frac{r_1}{a})f_{\alpha,\alpha}(x)+\frac{r_1}{a}f_{\alpha,\alpha+1}(x),\quad x>0.
\end{equation}
Integrating the previous equation, one obtains $\phi_-(0)=(\mathbb E I +\frac{r_1}{a})A_{\alpha,\alpha}+\frac{r_1}{a}A_{\alpha,\alpha+1}$
where $\displaystyle A_{\alpha,\beta}:=\int_{0}^{+\infty}f_{\alpha,\beta}(x)dx$.
According to Relation (\ref{Phi_-(0)}), one has $\phi_-(0)=1-\frac{r_1}{\mu-\lambda}$, and hence
\begin{equation}\label{r_1,EI,1}
1-(\frac{1}{\mu-\lambda}+\frac{A_{\alpha,\alpha}+A_{\alpha,\alpha+1}}{a})r_1=A_{\alpha,\alpha} \mathbb E I.
\end{equation}
Set $\displaystyle B_{\alpha,\beta}:=\int_{0}^{+\infty}xf_{\alpha,\beta}(x)dx$; multiplying Equation ({\ref{v_-(x)_linear-exp}}) by $x$,
integrating it and remembering the fact that $\mathbb E I=\displaystyle \int_{0}^{+\infty}xv_-(x)dx$, one obtains
\begin{equation}\label{r_1,EI,2}
(1-B_{\alpha,\alpha})\mathbb E I=\frac{B_{\alpha,\alpha}+B_{\alpha,\alpha+1}}{a}r_1.
\end{equation}
Thanks to Relations (\ref{r_1,EI,1}) and (\ref{r_1,EI,2}), one can deduce the constants $r_1$ and $\mathbb E I$; using
Equation ({\ref{v_-(x)_linear-exp}}), one can then completely determine the function $v_-$.

Finally, using Equation (\ref{phi+s}), we obtain
\begin{equation}\label{v_+_lin_exp}
v_+(x)=r_1e^{-(\mu - \lambda)x},\quad x>0.
\end{equation}
It is not surprising that the density of the workload, when positive, is exponentially distributed
with the same rate $\mu - \lambda$ as in the corresponding $M/M/1$ queue (arrival rate $\lambda$, service requirements exp($\mu$))
{\em without} inventory.
Indeed, every time the workload becomes positive, this occurs because of a customer arrival, and the memoryless property implies that the
residual part of the service requirement which makes that workload positive is $\exp(\mu)$ distributed.\\
\section {Direct approach}
\label{direct}
In this section, we use the analysis developed in \cite{AlbrecherBoxmaKuijstermans} to state some explicit results
when $\omega(x)=ax$ in the exponential service requirements case.
In \cite{AlbrecherBoxmaKuijstermans}, the authors study directly the functions $v_+$ and $v_-$ without considering their Laplace transforms. Indeed, differentiating Equations (\ref{queuing model>0}) and (\ref{queuing model<0}), one can show that the functions $v_-$ and $v_+$ satisfy some well known differential equations.

Set $C=\displaystyle \int_0^{+\infty} e^{-\mu x}v_-(x)dx$.
Differentiating Equation (\ref{queuing model>0}), one has
\begin{equation}\label{v_+_expression}
v_+(x)=C\lambda e^{-(\mu-\lambda)x}, \quad \forall \ x\geq 0.
\end{equation}
This is in agreement with
(\ref{v_+_lin_exp}), which we obtained following a Laplace transform approach.
On the other hand, differentiating Equation (\ref{queuing model<0}), one gets the following equation for $v_-$:
\begin{equation}\label{Eq_diff_v-_1}
v_-'(x)+(\lambda+ax)v_-(x)-\lambda \mu e^{\mu x} \int_x^{+\infty}e^{-\mu x}v_-(x)dx=0.
\end{equation}
Now, differentiating the expression in  Equation (\ref{Eq_diff_v-_1}), the function $v_-$ satisfies the following second order differential equation
\begin{equation}\label{Eq_diff_v-}
v_-''(x)+(\lambda-\mu+ax)v_-'(x)+a(1-\mu x)v_-(x)=0.
\end{equation}
Introduce the function $\theta(x)=v_-(x)e^{\frac{ax^2}{2}+\lambda x}$. The function $v_-$ is a solution of
Equation (\ref{Eq_diff_v-}) if and only if the function $\theta$ is a solution of the following second order differential equation:
\begin{equation}\label{Eq_diff_theta-}
\theta''(x)-(\lambda+\mu+ax)\theta'(x)+ \lambda \mu\theta(x)=0.
\end{equation}
One can check that $\theta$ is a solution of Equation (\ref{Eq_diff_theta-}) if and only if
$\theta(x)= \mathcal J \left( \tilde a, \tilde b,\frac{a}{2}(x+\frac{\lambda+\mu}{a})^2\right),$  where
$\tilde a=-\frac{\lambda \mu}{2a}$ and $\tilde b=\frac{1}{2}$ and $\mathcal J(\tilde a, \tilde b, \cdot)$ is a solution
of the degenerate hypergeometric equation:
\begin{equation}\label{degenerate_hypergeo}
zy''(z)+(\tilde b-z)y'(z)-\tilde a y(z)=0.
\end{equation}
According to \cite{NIST} (page 322)
and \cite{Poz} (page 143), Equation (\ref{degenerate_hypergeo}) has two standard solutions denoted by $z \mapsto M(\tilde a,\tilde b, z)$
and $z \mapsto U(\tilde a,\tilde b, z)$, the so called Kummer functions. Provided that $\tilde b \notin \{-1,-2,\cdots\}$, the function $z \mapsto M(\tilde a,\tilde b, z)$ is given by
$$M(\tilde a,\tilde b, z)=\displaystyle\sum _{0}^{+\infty}\frac{(\tilde a)_s}{(\tilde b)_s s!}z^s \quad \text{for all} \ z \in \mathbb C;$$
where $(c)_s= c(c+1) \cdots (c+s-1)$. The function $U(\tilde a,\tilde b,z)$ is uniquely determined by the property
$U(\tilde a,\tilde b,z)\sim z^{-a}$ when $z$ goes to $+\infty$. In our case ($\tilde b=\frac{1}{2}$), one has
\begin{equation}\label{relation_M_U}
U(\tilde a,\tilde b,z)=\frac{\Gamma(\frac{1}{2})}{\gamma(\frac{1}{2}-\frac{\lambda \mu}{2a})}M(-\frac{\lambda \mu}{2a},\frac{1}{2},z)+\Gamma (-\frac{1}{2})z^{\frac{1}{2}}M(-\frac{\lambda \mu}{2a}+\frac{1}{2},\frac{3}{2},z).
\end{equation}
The following analysis is developed in detail in \cite{AlbrecherBoxmaKuijstermans}, using knowledge of the degenerate hypergeometric equation (see for example \cite{NIST} page 322). It appears that one needs to distinguish between two cases:\\
Case $1$: $\tilde a=-\frac{\lambda \mu}{2a} \notin \mathbb Z,$ i.e $\frac{\lambda \mu}{2a} \notin \mathbb N$.
Denote by $\nu:=\frac{\lambda \mu}{a}.$
In this case, according to \cite{{AlbrecherBoxmaKuijstermans}}, one has
 \begin{equation}\label{v_(x)_final}
v_-(x)=Ke^{-\frac{a}{2}x^2-\lambda x}U\left(-\frac{\nu}{2},\frac{1}{2},\frac{a}{2}(x+\frac{\lambda+\mu}{a})^2\right) \quad \forall \ x\geq 0,
\end{equation}
where $K \in \mathbb R$. Relations (\ref{v_+_expression}) and (\ref{v_(x)_final}) imply that to determine completely the functions $v_+$ and $v_-$,
it is enough to determine the constants $C$ and $K$.
According to \cite{AlbrecherBoxmaKuijstermans}, these constants are given by the following system of equations:
\begin{eqnarray}
 &C=\frac{K(\lambda+\mu)U'(-\frac{\nu}{2},1/2, \frac{(\lambda+\mu)^2}{2a})}
 {\lambda \mu}, \label{Relation_K_2_C_1} \\
 &C\frac{\rho}{1-\rho}+K\int_0^{+\infty}e^{-\frac{a}{2}x^2-\lambda x}U(-\frac{\nu}{2},1/2,a/2(x+\frac{\lambda+\mu}{a})^2)dx=1.
 \label{Relation_K_2_C_2}
 \end{eqnarray}
 In particular, if $\nu$ is odd, i.e. $\nu=2n+1$, where $n \in \mathbb N$, from \cite{AlbrecherBoxmaKuijstermans}, one has
 \begin{equation}\label{v_(x)_odd}
 v_-(x)=2^{-\frac{\nu}{2}}K e^{-\frac{a}{2}x^2-\lambda x}H_{\nu}\left(\sqrt a x+\frac{\lambda+\mu}{\sqrt a}\right), \quad \forall \ x\geq 0,
\end{equation}
where, for $m \in \mathbb Z$, $H_m(\cdot)$ is the Hermite polynomial of order $m$ (see \cite{AS}, Page 775) given by the formula
$H_m(x)=(-1)^me^{\frac{1}{2}x^2}(\frac{d}{dx})^m[e^{\frac{-1}{2}x^2}]$.\\
Case $2$: $\frac{\lambda \mu}{2a}=n \in \mathbb N,$ i.e. $\nu=2n$.
In this case, according to \cite{AlbrecherBoxmaKuijstermans}, one has, for all $x\geq 0$,
\begin{equation}\label{Final_v_-_Hermite}
 v_-(x)=(-1)^{\frac{\nu}{2}}\frac{(\frac{\nu}{2})!}{\nu !}2^{\frac{\nu}{2}} K^*e^{-\frac{a}{2}x^2-\lambda x}H_{\nu}\left(\sqrt a x+\frac{\lambda+\mu}{\sqrt a}\right),\quad x \geq 0,
\end{equation}
where $K^* \in \mathbb R$.
Furthermore, similarly to case $1$, one can find two linear equations involving the unknowns $C$ and $K^*$ and then determine them.
One can see that Equations (\ref{v_(x)_odd}) and (\ref{Final_v_-_Hermite}) have the same shape.
Notice that Equation (\ref{Final_v_-_Hermite}) coincides with Equation (\ref{A_68}), the latter one being obtained in the appendix by a different method.
\section{Appendix: Analysis of first order differential equation}
We consider the differential equation (\ref{Eq-diff_Phi_-_1_exp})
\begin{equation}{\label{A_1}}
a (\mu-z)\frac{d}{dz}[\phi_-(z)]+z(z+\lambda-\mu)\phi_-(z)+a \mathbb E I(\mu-z)-r_1z=0\quad \text{for}
\quad \rm{Re}\ z\geq 0,
\end{equation}
where $a, \lambda, \mu$ are known positive constants with $\mu > \lambda$, and $\mathbb E I$ and $r_1$ are unknown. It was argued in
Section $2.1$ that $\phi_-$ is the Laplace transform of a non-negative, integrable function $v_-(x)$, $x\geq 0$, with $\int_0^{+\infty}v_-(x)dx<1$.
Hence, $\phi_-(z)$ is continuous, positive for $z \geq 0$ and analytic in the half plane for ${\rm{Re}}\ z>0$. Furthermore, $\phi_-(z)\rightarrow 0$ as
$ z\rightarrow \infty$. The objective is to compute $\mathbb E I$ and $r_1$ and to determine $\phi_-(z)$, from this information and from the fact that
\begin{equation}{\label{A_2}}
-\mathbb EI=\phi_-(0)=1-\frac{r_1}{\mu-\lambda}.
\end{equation}
The first equality in ({\ref{A_2}}) follows directly from ({\ref{A_1}}) by setting $z=0$, and the second one follows from Equation
(\ref{Phi_-(0)}) (see Section~\ref{linearomega}).\\
\

The standard calculus solution,
\begin{equation}{\label{A_3}}
y(x)=be^{-A(x)}+e^{-A(x)}\int_a^x Q(t)e^{A(t)}dt
\end{equation}
of the problem
\begin{equation}{\label{A_4}}
y'(x)+P(x)y(x)=Q(x), \quad y(a)=b,
\end{equation}
with
\begin{equation}{\label{A_5}}
 A(x)=\int_a^xP(t)dt,
\end{equation}
is  awkward to use since Equation ({\ref{A_1}}) gives rise to functions $P$ and $Q$ that are non-integrable around $z=\mu$.
Below, we will use a dedicated form
of the standard solution, separately on the ranges $0\leq z\leq \mu$ and $ z \geq \mu$.
\subsection{Considerations for $0\leq z \leq \mu$}\label{sub_sec_A_1}
We let
\begin{equation}{\label{A_6}}
  w =\mu-z \in [0,\mu];\ \chi( w)=\phi_-(\mu- w)=\displaystyle\sum_{k=0}^{\infty} c_k  w^k,
 \end{equation}
where the series converges for $| w|<\mu$. It follows easily from the fact that $\phi_-$ is the Laplace transform of an integrable non-negative function that all $c_k>0$. With notation ({\ref{A_6}}), Equation ({\ref{A_1}}) becomes
 \begin{equation}{\label{A_7}}
 a w \chi'( w)-(\mu- w)(\lambda - w) \chi( w)=\alpha +\beta  w,
 \end{equation}
 where the auxiliary quantities $\alpha$ and $\beta$ are defined by
  \begin{equation}{\label{A_8}}
 \alpha=-r_1 \mu,\ \beta=a \mathbb EI +r_1.
 \end{equation}
 From Equation ({\ref{A_7}}), there follows the recursion
 \begin{eqnarray}
 k=0:&\ &-\mu\lambda c_0=\alpha, {\label{A_9}}\\
 k=1:&\ &(a-\mu\lambda) c_1+ (\mu+\lambda)c_0=\beta,{\label{A_10}}\\
 k=2,3,\cdots:&\ &(ka-\mu\lambda) c_k+ (\mu+\lambda)c_{k-1}-c_{k-2}=0. {\label{A_11}}
 \end{eqnarray}
 We set
 \begin{equation}{\label{A_12}}
  U( w)=\sum_{k=K}^{+\infty}c_kw^{k}=\chi( w)-\sum_{k=0}^{K-1}c_k w^k,
 \end{equation}
where $K=1,2,\cdots$ is such that
\begin{equation}{\label{A_13}}
 K-1<\sigma:=\frac{\mu \lambda}{a} \leq K,
  \end{equation}
i.e., $K=\lceil \mu \lambda /a \rceil =\lceil \sigma \rceil.$ From Equation ({\ref{A_7}}), we have that
\begin{equation}{\label{A_14}}
\alpha +\beta  w=-\left(\mu \lambda-(\mu+\lambda) w\right)\chi(0)-(1-\frac{1}{\sigma})\mu\lambda \chi'(0) w.
\end{equation}
Therefore, Equation ({\ref{A_7}}) can be rewritten as
\begin{equation}{\label{A_15}}
\frac{\mu \lambda}{\sigma}  w\left(\chi( w)-\chi(0)-\chi'(0) w\right)'-\mu \lambda \left(\chi( w)-\chi(0)-\chi'(0) w\right)
+(\mu+\lambda) w\left(\chi( w)-\chi(0)\right)- w^2\chi( w)=0.
\end{equation}
Writing the functions $\chi( w)$, $\chi( w)-\chi(0)$ and $\chi( w)-\chi(0)-\chi'(0) w$ that occur in
Equation ({\ref{A_15}}) in terms of $U( w)$ and $\displaystyle\sum_kc_k w^k$, a careful administration with the recursion in
Relations ({\ref{A_9}-{\ref{A_11}}) gives
\begin{equation}{\label{A_16}}
\frac{\mu \lambda}{\sigma} ( w U'( w)-\sigma U( w))+ ((\mu+\lambda) w- w^2)U( w)
+\frac{\sigma-K}{\sigma}\mu\lambda c_K w ^K-c_{K-1} w^{K+1}=0.
 \end{equation}
Observing that $ w U'( w)-\sigma U( w)=w^{\sigma+1}( w^{-\sigma}U( w))'$ and dividing Equation ({\ref{A_16}})
 by $(\mu\lambda/\sigma) w^{\sigma+1}$, we get
 \begin{equation}{\label{A_17}}
 (\frac{U( w)}{ w^{\sigma}})' +\frac{\sigma(\mu+\lambda- w)}{\mu \lambda}\frac{U( w)}{ w^{\sigma}}
 =-(\sigma-K)c_K w^{K-\sigma-1}+\frac{\sigma}{\mu\lambda}c_{K-1} w^{K-\sigma}.
 \end{equation}
 The right hand side of  Equation ({\ref{A_17}}) is integrable at $ w=0$ because of the choice of $K$ in Equation ({\ref{A_13}}).
 Hence, the formula ({\ref{A_3}}) can be used with $y( w)= w^{-\sigma}U( w)$ and $a=0$, where we note that
  \begin{eqnarray}{\label{A_18}}
  \frac{U( w)}{ w^{\sigma}}|_{ w=0}=
  \left\{
    \begin{split}
    &0\quad,\quad K-1<\sigma<K,\\
    &c_K\quad,\quad \sigma=K.
    \end{split}
  \right.
\end{eqnarray}
This yields for $K-1<\sigma<K$:
\begin{equation} \label{A_19}
\frac{1}{ w^{\sigma}}U( w)=e^{-A_1( w)}\int_0^{ w}Q_1(v)e^{A_1(v)}dv,
\end{equation}
and for $\sigma=K$:
\begin{equation} \label{A_20}
\frac{1}{ w^{K}}U( w)=c_Ke^{-A_1( w)}+e^{-A_1( w)}\int_0^{ w}Q_1(v)e^{A_1(v)}dv,
\end{equation}
 where
 \begin{equation}\label{A_21}
 A_1( w)=\int_0^{ w}\frac{\sigma(\mu+\lambda-v)}{\mu\lambda}dv=\frac{\sigma}{\mu\lambda}((\mu+\lambda) w-\frac{1}{2} w^2),
 \end{equation}
and
\begin{equation}\label{A_22}
Q_1(v)=-(\sigma-K)c_Kv^{K-\sigma-1}+\frac{\sigma}{\mu\lambda}c_{K-1}v^{K-\sigma},
\end{equation}
with a constant $Q_1(v)=\frac{\sigma c_{K-1}}{\mu\lambda}$ when $\sigma=K$. Also, observe
from Equation (\ref{A_13}) that $\frac{\sigma}{\lambda \mu}=a$ in Equations (\ref{A_21}-\ref{A_22}).
\subsection{Considerations for $z\geq \mu$}\label{sub_sec_A_2}
We now let
\begin{equation}\label{A_23}
 w=z-\mu \geq 0;\quad \Lambda( w)=\phi_-(\mu+ w)=\sum_{k=0}^{\infty}d_k w^k,
\end{equation}
where $d_k=(-1)^kc_k,$ and we set
\begin{equation}\label{A_24}
 T( w)=\sum_{k=K}^{+\infty}d_k w^k=\Lambda( w)-\sum_{k=0}^{K-1}d_k w^k.
\end{equation}
With a similar analysis as in Subsection \ref{sub_sec_A_1}, we get for $K-1<\sigma<K$:
\begin{equation}\label{A_25}
\frac{1}{ w^K}{T( w)}=e^{-A_2( w)}\int_0^{ w}Q_2(v)e^{A_2(v)}dv,
\end{equation}
and for $\sigma=K$:
\begin{equation}\label{A_26}
\frac{1}{ w^K}{T( w)}=d_Ke^{-A_2( w)}+e^{-A_2( w)}\int_0^{ w}Q_2(v)e^{A_2(v)}dv,
\end{equation}
where
\begin{equation}\label{A_27}
 A_2( w)=-\int_0^{ w} \frac{\sigma(\mu+\lambda +v)}{\mu \lambda}dv=-\frac{\sigma}{\mu\lambda}\left((\mu+\lambda) w+\frac{ w^2}{2}\right).
\end{equation}
and
\begin{equation}\label{A_28}
  Q_2(v)=-(\sigma-K)d_Kv^{K-\sigma-1}+\frac{\sigma}{\mu\lambda}d_{K-1}v^{K-\sigma},
\end{equation}
with a constant $Q_2(v)=\frac{\sigma d_{K-1}}{\mu\lambda}$ when $\sigma=K.$
\subsection{Boundary condition at $z=0$}\label{sub_sec_A_3}
We have from
\begin{equation}\label{A_29}
 \phi_-(\mu- w)=\chi( w)=U( w)+\sum_{k=0}^{K-1}c_k w^k,
\end{equation}
and Equations (\ref{A_19}-\ref{A_20}) that
\begin{equation}\label{A_30}
 \phi_-(\mu- w)=\sum_{k=0}^{K-1}c_k w^k+c_K w^Ke^{-A_1( w)}+ w^\sigma e^{-A_1( w)}
 \int_0^ w Q_1(v)e^{A_1(v)}dv,
\end{equation}
where the term comprising $c_K$ is to be included only in the case that $\sigma=K$. The left-hand side
of Equation (\ref{A_30}) equals $ 1-\frac{r_1}{\mu-\lambda}$ at $ w=\mu$ by Equation (\ref{A_2}), and so we get
\begin{equation}\label{A_31}
  1-\frac{r_1}{\mu-\lambda}=\sum_{k=0}^{K-1}c_k\mu^k+c_K\mu^Ke^{-A_1(\mu)}+\mu^\sigma e^{-A_1(\mu)}\int_0^\mu Q_1(v)e^{A_1(v)}dv,
\end{equation}
with the term comprising $c_K$ to be included only when $\sigma=K$.
\subsection{Boundary condition at $z=+\infty$}\label{sub_sec_A_4}
We have from Equation (\ref{A_24}) and the fact that $\phi_-(z)\rightarrow 0$ as $z\rightarrow +\infty$ that
$T( w)=O( w^{K-1})$ as $ w \rightarrow +\infty$. Therefore, for $K-1<\sigma<K$, we have that
\begin{equation}\label{A_32}
\int_0^{+\infty}Q_2(v)e^{A_2(v)}dv=0,
\end{equation}
since the left-hand side of Equation (\ref{A_25}) goes to $0$ as $ w \rightarrow +\infty$. Using Equation
(\ref{A_28}) and $d_K=(-1)^Kc_{K}$, $d_{K-1}=(-1)^{K-1}c_{K-1},$ this gives
\begin{equation}\label{A_33}
 c_{K-1}=\frac{\mu\lambda}{\sigma}(K-\sigma)\frac{I_{K-\sigma-1}}{I_{K-\sigma}}c_K,
\end{equation}
where
\begin{equation}\label{A_34}
 I_{L-\sigma}=\int_0^{+\infty}v^{L-\sigma}e^{A_2(v)}dv, \quad L=K,K-1.
\end{equation}
 Similarly, when $\sigma=K$ (noting the simplification of $Q_2$ in Equation (\ref{A_28}))
 \begin{equation}\label{A_35}
 c_{k-1}=\frac{\mu\lambda}{\sigma}\frac{1}{I_0}c_K.
\end{equation}
\subsection{Rewriting all unknowns in terms of $c_K$}\label{sub_sec_A_5}
We have from Equations (\ref{A_8}) and (\ref{A_9}-\ref{A_10})
\begin{equation}\label{A_36}
\begin{bmatrix}
r_1 \\
\mathbb E I
\end{bmatrix}=
\begin{bmatrix}
-\frac{1}{\mu} & 0 \\
\frac{1}{a\mu} & \frac{1}{a}
\end{bmatrix}
\begin{bmatrix}
 \alpha \\
\beta
\end{bmatrix},
\quad
\begin{bmatrix}
 \alpha \\
\beta
\end{bmatrix}
=
\begin{bmatrix}
-\mu\lambda & 0 \\
 \mu+\lambda & a-\mu\lambda
\end{bmatrix}
\begin{bmatrix}
 c_0\\
 c_1
\end{bmatrix},
\end{equation}
and so
\begin{equation}\label{A_37}
\begin{bmatrix}
r_1 \\
\mathbb E I
\end{bmatrix}=
\begin{bmatrix}
 \lambda & 0 \\
 \frac{\mu}{a} &  1-\frac{\mu\lambda}{a}
\end{bmatrix}
\begin{bmatrix}
 c_0 \\
 c_1
\end{bmatrix}.
\end{equation}
Furthermore, from the first item in Equations (\ref{A_8}) and (\ref{A_9}),
\begin{equation}\label{A_38}
1-\frac{r_1}{\mu-\lambda}=1+\frac{\alpha/\mu}{\mu-\lambda}=1-\frac{\lambda}{\mu-\lambda}c_0.
\end{equation}
We next use the recursion in Equation (\ref{A_11}) in backward order to write all $c_k,\ k=K,K-1,\cdots,0$
in terms of $c_K,$ where for the case that $k=K-1$, we have by equations (\ref{A_33}), (\ref{A_35}),
\begin{eqnarray}\label{A_39}
c_{K-1}=B_{K-1}c_K;\quad B_{K-1}=
  \left\{
    \begin{split}
    &\frac{\mu\lambda}{\sigma}(K-\sigma)\frac{I_{K-\sigma-1}}{I_{K-\sigma}},\quad K-1<\sigma<K,\\
    &\frac{\mu\lambda}{\sigma}\frac{1}{I_0},\quad \sigma=K.
    \end{split}
  \right.
\end{eqnarray}
We find
\begin{equation}\label{A_40}
 c_{k-2}=B_{k-2}c_K,\ k=K,K-1,\cdots,2,
\end{equation}
where the $B_k$ are given recursively according to
\begin{equation}\label{A_41}
 B_{k-2}=(ka-\mu\lambda)B_k+(\mu+\lambda)B_{k-1},\ k=K,K-1,\cdots,2,
\end{equation}
with $B_K=1$ and $B_{K-1}$ given in Equation (\ref{A_39}).
\subsection{Solving for $c_K$ and finding $r_1$ and $\mathbb E I$}\label{sub_sec_A_6}\ \\
We have from Equations (\ref{A_21}-\ref{A_22})
\begin{equation}\label{A_42}
 \int_0^\mu Q_1(v)e^{A_1(v)}dv=-(\sigma-K)\mathcal J_{K-\sigma-1}c_K+\frac{\sigma}{\mu\lambda}\mathcal J_{K-\sigma}c_{K-1},
\end{equation}
where the term with $c_K$ vanishes when $\sigma=K$, and
\begin{equation}\label{A_43}
 \mathcal J_{L-\sigma}=\int_0^\mu v^{L-\sigma}e^{A_1(v)}dv,\quad L=K,K-1.
\end{equation}
The right-hand side of Equation (\ref{A_31}) can then be written as
\begin{equation}\label{A_44}
 \sum_{k=0}^{K-1}c_k\mu^k+\mu^{\sigma}e^{-A_1(\mu)}\frac{\sigma}{\mu\lambda}\mathcal J_{K-\sigma}c_{K-1}
 +\mu^{\sigma}e^{-A_1(\mu)}c_KR_K,
\end{equation}
where
\begin{eqnarray}\label{A_45}
 R_K=
  \left\{
    \begin{split}
     &-(\sigma-K)\mathcal J_{K-\sigma-1},\quad K-1<\sigma<K,\\
     &1,\quad\sigma=K.
    \end{split}
  \right.
\end{eqnarray}
For the left-hand side of Equation (\ref{A_31}) we use Relation (\ref{A_38}) and we bring the term
$\frac{\lambda c_0}{(\mu-\lambda)}$ to the right-hand side. Then using $c_k=B_kc_K$ and treating the term
with $k=0$ separately, we find for Relation (\ref{A_31})
\begin{equation}\label{A_46}
1=\left(\frac{\mu}{\mu-\lambda}B_0+\sum_{k=1}^{K-1}B_k\mu^k+
\mu^\sigma e^{-A_1(\mu)}[\frac{\sigma}{\mu\lambda}\mathcal J_{K-\sigma}B_{K-1}+R_K]\right)c_K.
\end{equation}
In Equation (\ref{A_46}), all quantities in the parentheses are known, and so we have found $c_K$. Then, with Equations (\ref{A_39}) and
(\ref{A_40}) all $c_k,\ k=0,1,\cdots,K-1,$ can be found, and, finally, Equation (\ref{A_37}) yields $r_1$ and $\mathbb EI$.

\subsection{Special attention for the case $0<\sigma\leq 1$}\label{sub_sec_A_7}\ \\
When $0<\sigma<1$, we can find $c_1$ and $c_0$ from Equations (\ref{A_46}) and (\ref{A_39}). We note here that
\begin{equation}\label{A_47}
 B_0=\frac{\mu\lambda}{\sigma}(1-\sigma)\frac{I_{-\sigma}}{I_{1-\sigma}},\quad R_1=(1-\sigma){\mathcal J}_{-\sigma}.
\end{equation}
With $\frac{\mu\lambda}{\sigma}=a$, see Equation (\ref{A_13}), we then get
\begin{equation}\label{A_48}
 c_1=\frac{\frac{1}{1-\sigma}I_{1-\sigma}}{\frac{\mu a}{\mu-\lambda}I_{-\sigma}+\mu^{\sigma}e^{-A_1(\mu)}(\mathcal J_{1-\sigma}
 I_{-\sigma}+I_{1-\sigma}\mathcal J_{-\sigma})},
\end{equation}
\begin{equation}\label{A_49}
 c_0=B_0c_1=\frac{aI_{-\sigma}}{\frac{\mu a}{\mu-\lambda}I_{-\sigma}+\mu^{\sigma}e^{-A_1(\mu)}(\mathcal J_{1-\sigma}
 I_{-\sigma}+I_{1-\sigma}\mathcal J_{-\sigma})}.
\end{equation}
Finally, $r_1=\lambda c_0,$ and $\mathbb EI$ is obtained from Relation (\ref{A_37}) as
\begin{equation}\label{A_50}
 \mathbb E(I)=\frac{\mu}{a}c_0+(1-\sigma)c_1=\frac{\mu I_{-\sigma}+I_{1-\sigma}}{\frac{\mu a}{\mu-\lambda}I_{-\sigma}+\mu^{\sigma}e^{-A_1(\mu)}(\mathcal J_{1-\sigma}
 I_{-\sigma}+I_{1-\sigma}\mathcal J_{-\sigma})}.
\end{equation}
From $0<\lambda<\mu$ and
\begin{equation}\label{A_51}
\frac{r_1}{\lambda}=c_0,\quad \frac{r_1}{\mu-\lambda}=\frac{r_1\lambda}{\lambda (\mu-\lambda)},
\end{equation}
we see that Equation (\ref{A_49}) yields that both quantities in Relation (\ref{A_51}) are between $0$ and $1$ (see Equation (\ref{A_2})).

In the case that $\sigma=1$, we have that $B_0=\frac{a}{I_0},\ R_0=1$,
and we get consequently
\begin{equation}\label{A_52}
 c_1=\frac{I_0}{\frac{\mu a}{\mu-\lambda}+\mu e^{-A_1(\mu)}(\mathcal J_0+I_0)},
\end{equation}
\begin{equation}\label{A_53}
 c_0=B_0c_1=\frac{a}{\frac{\mu a}{\mu-\lambda}+\mu e^{-A_1(\mu)}(\mathcal J_0+I_0)}=\frac{r_1}{\lambda},
\end{equation}
\begin{equation}\label{A_54}
\mathbb E(I)=\frac{\mu c_0}{a} = \frac{1}{\frac{a}{\mu-\lambda}+  e^{-A_1(\mu)}(\mathcal J_0+I_0)}.
\end{equation}
\subsection{Laplace transform representation of $\phi_-$ when $\sigma=K$}\label{sub_sec_A_8}\ \\
In the case that $\sigma=K=1,2,\cdots$, we have from Subsection \ref{sub_sec_A_2} that
\begin{equation}\label{A_55}
 \phi_-(\mu+ w)=T( w)+\sum_{k=0}^{K-1}d_k w^k= w^Ke^{-A_2( w)}\left[d_K+\frac{K}{\mu\lambda}d_{K-1}\int_0^{ w}e^{A_2(v)}dv\right]
 +\sum_{k=0}^{K-1}d_k w^k.
\end{equation}
From Equation (\ref{A_35}) and $d_K=(-1)^Kc_K,\ d_{K-1}=(-1)^{K-1}c_{K-1}$, we find
\begin{equation}\label{A_56}
 d_K=-\frac{K}{\mu\lambda}d_{K-1}\int_0^{+\infty}e^{A_2(v)}dv,
 \end{equation}
and this yields
\begin{equation}\label{A_57}
\phi_-(\mu+ w)=-\frac{K}{\mu\lambda}d_{K-1} w^K e^{-A_2( w)}\int_ w^{+\infty}e^{A_2(v)} dv+\sum_{k=0}^{K-1}d_k w^k.
\end{equation}
By analyticity of $\phi_-(z)$ in $z>0$ and the right-hand side in $ w \in \mathbb C$, we have that Equation (\ref{A_57})
holds for all $ w \geq -\mu$. We manipulate now the right-hand side of Equation (\ref{A_57}) to find $\phi_-(z)$
as the Laplace transform of a non-negative function $f(t)$ -- which of course is $v_-(t)$.
With $z=\mu+ w\geq 0$ and $v=u+ w$, we have from Equation (\ref{A_27})
\begin{equation}\label{A_58}
 -A_2( w)+A_2(v)=\frac{K}{\mu\lambda}\left[(\mu+\lambda) w+\frac{1}{2} w^2-(\mu+\lambda)v-\frac{1}{2}v^2\right]=\frac{-K}{2\mu\lambda}
\left[2zu+2\lambda u+u^2\right].
\end{equation}
Therefore,
\begin{equation}\label{A_59}
 e^{-A_2( w)}\int_{ w}^{+\infty}e^{A_2(v)}dv=\int_0^{+\infty}\exp(-\frac{K}{2\mu\lambda}[2zu+2\lambda u+u^2])du
 =\frac{\mu\lambda}{K}\int_0^{+\infty} e^{-zt}\exp(-\lambda t -\frac{\mu \lambda}{2K}t^2)dt,
\end{equation}
where for the second step the substitution $u=\frac{\mu \lambda t}{K}$ has been made. Setting
\begin{equation}\label{A_60}
 \Psi_K(t)=\exp(-\lambda t -\frac{\mu \lambda}{2K}t^2),
\end{equation}
we thus have for $z=\mu+ w \geq 0$
\begin{equation}\label{A_61}
 \phi_-(z)= -d_{K-1}(z-\mu)^K\int_0^{+\infty}e^{-zt}\Psi_K(t)dt+\sum_{k=0}^{K-1}d_k(z-\mu)^k.
\end{equation}
From basic properties of the Laplace transform, we have
\begin{equation}\label{A_62}
 z^j\int_{0}^{+\infty}e^{-zt}\Psi_K(t)dt=\sum_{i=0}^{j-1}z^{j-i-1}\Psi_K^{(i)}(0)+\int_{0}^{+\infty}e^{-zt}\Psi_K^{(j)}(t)dt,
\end{equation}
for $j=0,1,\cdots .$ Hence, by Newton's binomium for $(z-\mu)^K,$ we get
\begin{equation}\label{A_63}
 \phi_-(z)= -d_{K-1} \int_0^{+\infty}e^{-zt} \sum_{j=0}^{K}\binom{K}{j}(-\mu)^{K-j}\Psi_K^{(j)}(t)dt+P(z),
\end{equation}
 where $P(z)$ is a polynomial. Since $\phi_-(z)$ and the integral at the right-hand side of Equation (\ref{A_63})
go to $0$ as $z\rightarrow +\infty,$ we have $P(z)=0$. It is a consequence of the formula
\begin{equation}\label{A_64}
\sum_{j=0}^{K}\binom{K}{j}(-\mu)^{K-j}\Psi_K^{(j)}(t) = (\frac{d}{ds})^K[e^{-\mu s}\Psi_K(t+s)]_{s=0}
\end{equation}
and the generating function
\begin{equation}\label{A_65}
\sum_{n=0}^{+\infty}H_n(x)\frac{z^n}{n!} = e^{xz-\frac{1}{2}z^2}
\end{equation}
of the Hermite polynomials $H_n(x)=(-1)^ne^{\frac{1}{2}x^2}(\frac{d}{dx})^n[e^{\frac{-1}{2}x^2}]$
that
\begin{equation}\label{A_66}
 \sum_{j=0}^{K}\binom{K}{j}(-\mu)^{K-j}\Psi_K^{(j)}(t)= (-1)^K (\frac{\mu \lambda}{K})^{\frac{K}{2}}\Psi_K(t) H_K\left((\frac{K}{\mu \lambda})^{\frac{1}{2}}(\mu+\lambda)+(\frac{\mu \lambda}{K})^{\frac{1}{2}} t \right).
\end{equation}
Then using that $d_{K-1}=(-1)^{K-1}c_{K-1},$ we finally obtain
\begin{equation}\label{A_67}
\phi_-(z)=\int_0^{+\infty}e^{-zt}f(t)dt,
\end{equation}
where
\begin{equation}\label{A_68}
 f(t)= c_{K-1}(\frac{\mu \lambda}{K})^{\frac{K}{2}}\Psi_K(t) H_K\left((\frac{K}{\mu \lambda})^{\frac{1}{2}}(\mu+\lambda)+(\frac{\mu \lambda}{K})^{\frac{1}{2}} t \right).
\end{equation}
This displays $\phi_-$ in the form of a Laplace transform, and  we have found $v_-$ which equals $f$. One can notice that Equation
(\ref{A_68}) agrees with the result obtained previously in Equation (\ref{Final_v_-_Hermite}).\\
From
\begin{equation}\label{A_69}
 c_{K-1}>0,\ \frac{\mu+\lambda}{(\mu \lambda)^{\frac{1}{2}}}\geq 2;\ H_K(x)>0,\ x>2\sqrt K,
\end{equation}
we see that $f(t)$ in Equation (\ref{A_68}) is positive for $t\geq 0$. Here we have used \cite{Zeg}, Theorem 6.32 on pp $131-132$ that implies
that the largest zero of $H_n(x)$ is less than $\sqrt{2n+1}$.
\ \\

\noindent
{\bf Acknowledgment}
\\
The authors gratefully acknowledge stimulating discussions with Hansjoerg Albrecher.
The research of Onno Boxma was supported by the NETWORKS project, which is being funded by the Dutch government. The research of Rim Essifi
was supported by the ERC CritiQueue program.

\end{document}